\documentclass{amsart}
\usepackage{amsmath,amscd,amstext,amsthm,amsfonts,latexsym}
\usepackage{eufrak, graphicx,epstopdf,stmaryrd}
\usepackage{mathtools, abraces,yhmath, float, amssymb}
\usepackage{qtree, forest, tikz-qtree}
\usepackage{color, fancybox, empheq}
\usetikzlibrary{fit}
\usepackage{caption}

\DeclareCaptionType{InfoBox}
\usetikzlibrary{shapes}
\useforestlibrary{linguistics}
\usepackage{multirow, longtable, caption}

\theoremstyle{plain}
\newtheorem{theorem}{Theorem}[section]
\newtheorem{proposition}[theorem]{Proposition}
\newtheorem{lemma}[theorem]{Lemma}

\newtheorem{corollary}[theorem]{Corollary}
\newtheorem{maintheorem}{Theorem}

\newtheorem{remark}{Remark}

\newtheorem{definition}[theorem]{Definition}

\newtheorem{notation}{Notation}

\newcommand{\field}[1]{\mathbb{#1}}

\newcommand{\ZZ}{\field{Z}}
\newcommand{\NN}{\field{N}}

\newcommand{\FF}{\field{F}}
\def\N{\NN}
\def\ninf{n\to+\8}

\def\kinf{k\to+\8}

\newcommand{\al} {\alpha}       
\newcommand{\be} {\beta}        
\newcommand{\ga} {\gamma}

\newcommand{\la} {\lambda}

\newcommand{\om} {\omega}

\newcommand{\cO}{{\mathcal O}}

\def\s{\sigma}

\def\disp{\displaystyle}
\def\FFF{\FF_{2}^{+}}
\def\ie{{\em i.e.,\ }}

\def\8{\infty}
\def\wt{\widetilde}

\newcommand{\ol}[1]{\overline{#1}}

\newcommand {\CA}{{\mathcal A}}
\newcommand {\CB}{{\mathcal B}}

\newcommand {\CO}{{\mathcal O}}
\newcommand {\CP}{{\mathcal P}}

\newcommand {\CT}{{\mathcal T}}

\newcommand {\GA}{{\mathfrak A}}
\newcommand {\GB}{{\mathfrak B}}
\newcommand {\GC}{{\mathfrak C}}
\newcommand {\GD}{{\mathfrak D}}
\newcommand {\GE}{{\mathfrak E}}
\newcommand {\GF}{{\mathfrak F}}
\newcommand {\GG}{{\mathfrak G}}
\newcommand {\GH}{{\mathfrak H}}
\newcommand {\GI}{{\mathfrak I}}
\newcommand {\GJ}{{\mathfrak J}}
\newcommand {\GK}{{\mathfrak K}}
\newcommand {\GL}{{\mathfrak L}}
\newcommand {\GM}{{\mathfrak M}}
\newcommand {\GN}{{\mathfrak N}}

\newcommand {\GR}{{\mathfrak R}}

\newcommand{\eg}{\vskip 0.2cm \noindent{\bf Example:} \hskip 0.2cm}

\def\ST{\CA^{\FF_2^+}}
\newcommand{\indiceGauche}[2]{{\vphantom{#2}_{#1}#2}}
\newcommand{\ig}{\indiceGauche}
\def\afz{{\ig{0}{\GA}}}
\def\afu{{\ig{1}{\GA}}}
\def\afi{{\ig{i}{\GA}}}

\newcommand{\comment}[1]{}

\newcommand\RL{\color{red}}

\newcommand{\arb}[1]{{\begin{forest}
[$#1$[{}][{}]]
\end{forest}}}

\newcommand{\tritree}[3]{
\begin{forest}
[{#1}[{#2}][{#3}]]
\end{forest}
}

\newcommand{\heptatree}[7]{
\begin{forest}
[{#1}
[{#2}[{#4}][{#5}]]
[{#3}[{#6}][{#7}]]
]
\end{forest}
}

\def\bw{\textrm{bw}}

\def\vbw{\textrm{vbw}}

\begin{document}

\title{Substreetutions and more on trees}
\author{A. Baraviera}
\thanks{AB wants to thank ISEA for kind support for a visit at Noum\'ea in January 2019}
\address{Instituto de Matemática e Estatística - UFRGS\\
Avenida Bento Gonçalves, 9500  - Porto Alegre - RS - Brasil\\
CEP 91509-900}
\email{baravi@mat.ufrgs.br}

\author{Renaud Leplaideur}
\address{ISEA, Universit\'e de la Nouvelle-Cal\'edonie \& LMBA UMR6205}
\email{renaud.leplaideur@unc.nc}

\date{\today}


\maketitle

 \begin{abstract}
    We define a notion of substitution on colored binary trees that we call substreetution. We show that a fixed point by a substreetution may be (or not) almost periodic, thus the closure of the orbit under $\FFF$-action may (or not) be minimal.

    We study one special example: we show that it belongs to the minimal  case and  that the number of preimages in the minimal set increases just exponentially fast, whereas it could be expected a super-exponential growth.

    We also give examples of periodic trees without invariant measure on their orbit.

    We use our construction to get quasi-periodic colored tilings of the hyperbolic disk.
 \end{abstract}

 \tableofcontents

\section{Introduction}
\subsection{Background and main motivations}
In the present paper we define and study a notion of substitution on colored trees. In the following they will be called  substreetutions to make clear they are substitutions on trees. Substitution will be used for ``classical ones''  acting on the line. These substreetutions are new objects and generate many natural questions, within the dynamical viewpoint.

Our first main theorem (see Theorem \ref{MainThA}) states the existence of a substreetution $H$ such that the $\FFF$-orbit of the fixed point  (for $H$), $\GJ$, generates a minimal set $X$ in the set of colored binary trees (under $\FFF$ action).

\bigskip
The first motivation comes from Thermodynamic formalism in Ergodic Theory. A series a works (\cite{BLL,BL1,BL2}) investigates the question of freezing phase transition with a quasi-crystal ground state and make links with substitutions. The question of higher-dimensional actions naturally appears, but conflicts to the fact that 1-d thermodynamical formalism deeply uses the Transfer Operator. This operator needs (to be defined) a kind of canonical order in the action, which does not hold for  $\ZZ^{d}$ actions with $d\ge 2$. Considering $\FFF$-actions was thus the way to recover such canonical order in the action.  We have thus been led to study substitutions associated to $\FFF$-actions.

\bigskip
The second motivation is the existence of works on Sturmian trees (see \cite{KLLS}  and \cite{BBCF}).
 It was a motivation to define substitutions on colored trees and to study the closure of the fixed point, as it is done for the Fibonacci sequence. We remind that the Fibonacci sequence is both Sturmian and substitutive.

\bigskip
The third motivation deals with the notion of invariant measures. Coming from the Ergodic Theory world, this question is of prime interest.
It seems that few is known on  invariant measures for trees. Unless we consider simple periodic trees with high symmetries, Bernoulli measures on the vertices or Markov inspired measures (see \cite{BP}, \cite{Spitzer}),  
we are not aware of other examples. We remind that $\FFF$ is not amenable, which makes the standard construction of invariant measures fail.
In \cite{FW}, some stationary measures are studied but it is not clear that they are $\FFF$-invariant measures. We also mention the work of Rozikov \emph{et al} (see \cite{Rozikov}) in the statistical mechanics context.

We remind that in the classical Thermodynamic Formalism,  invariant measures are usually obtained as \emph{equilibrium states},  and more precisely as \emph{conformal} measures coming from a Transfer Operator. Note that this operator exists here,  but the question remains to see if we can use it to construct conformal measures. In particular, does it act on continuous functions ?

If one wants to mimic or adapt the classical proofs, one of the first key points is to control the spectral radius of the Transfer Operator. Still for the classical case, the logarithm of the spectral radius is the  topological entropy of the system. Entropy is a way to quantify the complexity of the system, and this is  the core of the study in Dynamical Systems. Usually, entropy is given by the exponential growth rate of the number of configurations  the system sees with respect to their length. For trees, it is naturally expected a super-exponential growth rate (see \cite{Petersen-Salama-20}). In  \cite{KLLS}, Sturmian trees are colored trees with an affine growth.

At that point, two questions arise: what is the good renormalization scale to classify chaos, and how to compute entropy? We remind that the spectral radius of the Transfer Operator usually deals with backward iterates, whereas the number of
configurations/patches usually deals with forward iterates.

For the first question, we remind that for classical substitutions, complexity is in  $O(n)$, $O(n\log\log n)$, $O(n\log n)$ or $O(n^{2})$. This means that chaos is not observed at exponential scale but at a lower scale. Substitutions are chaotic but zero-entropy systems.
Our second main result Theorem \ref{MainThB} exhibites a similar result but with a change of scale. Whereas we could expect a super exponential growth rate, we show here that we \emph{only} have an exponential growth rate.

The second question has been studied for many systems. The result is that \emph{in spirit} both rates do coincide.
In particular, we mention \cite{Przytycki-pressure}, where several notion of pressures\footnote{We remind that the pressure is a kind of extension of the entropy.} are defined for rational quadratic maps. Some are computed involving backward iterates, some with forward iterates. In \cite{PRLS} it was proved that all these pressures do coincide.

In the present paper, we define a notion of entropy, inspired by the tree-pressure\footnote{Here tree-pressure is related to the tree-structure of inverse branches under a $\ZZ$-action and has not to be mixed up with the tree structure coming from $\FFF$-action.} from \cite{Przytycki-pressure}. We show in our second main result Theorem \ref{MainThB} that the growth rate for the number of preimages of any tree in $X$ is at most exponential, with a bounded from above speed. To avoid confusion with what would be entropy for trees, we decided to call it the \emph{backward complexity}. However, we conjecture that this backward complexity has the same growth than the number of patches that can be seen in $\GJ$.

\bigskip
 Copying the vocabulary for classical substitution, we \emph{mostly} consider  in the present paper constant length-2 substitutions on regular binary trees. More general substreetutions can easily be defined, however, some computations show that complexity may become extremely quickly very high.

Several examples of constant length-2 substreetutions are described and studied in the paper. For one of them,  the fixed point is the usual Thue-Morse sequence lifted in a tree (each line has a constant color). The closure of the orbit of this tree is obviously minimal but does not really exploit the tree structure.

Furthermore, in Theorem \ref{th-perio-nomeas}, we give one example of a periodic tree without invariant measure.

\subsection{Settings and results}

\subsubsection{ Binary trees and colored trees}
We consider the free monoid with two generators $\FF_2^+$ whose elements are the identity  and all the finite words, of any given length,
that can be written with the alphabet $\{ a, b \}$. The letter $e$ stands for the empty word.

 $| \cdot |$ is the length of a
given word, meaning the number of letters it has; we also set $|e|=0$.
 The product
of two words is just the word obtained by the concatenation. The lexicographic order on $\FFF$ is the one given by $a<b$ and $\omega<\om'$ if $|\om|<|\om'|$.

\bigskip
Fix a  set $\CA$, called the alphabet;
a {\it colored  tree (also called configuration)} is a map $\mathfrak B \colon \FF_2^+ \to \CA$, say $\mathfrak B \in \CA^{\FF_2^+}$. In the following one shall set $\mathfrak B_{x}$ for the value $\GB(x)$ and,
for simplicity, we also fix $\CA$ as being the set $\{0, 1\}$. We will only consider colored trees and thus simply use the word tree.

If $\GA$ is a tree and $\om$ belongs to $\FFF$, $\GA_{\om}$ is called digit at position $\om$. Hence, if the tree $\GA$ is fixed, elements of $\FFF$ will also be called {\it sites}. The {\it generation} of a site $\om$ is $|\om|$.
Hence, the {\it root}, that corresponds to the site $e$  is generation 0, the two first descendants form the set of sites at  generation 1 and so on.

If $\om$ is a site, $\om.a$ is called its $a$-follower. Equivalently, we shall say that the digit $A_{\om a}$ is  the $a$-follower of digit $A_{\om}$ and the digit  $A_{\om b}$ is  the $b$-follower of digit $A_{\om}$. The sites $\om.a$ and $\om.b$ are said to be brothers. Followers will also be called children. A $n$-descendant for $\om$ is any digit at site $\om\om'$ with $|\om'|=n$. Children are 1-descendants.

The site $\om$ is the \emph{1-ascendant} of sites $\om.a$ and $\om.b$. We shall also call it the \emph{father}.
The father of the father is called the \emph{grandfather}. It is also the 2-ancestor. By induction, one may define the $n$-ancestor of any site $\om$ with $|\om|>n$. Two sites are \emph{cousin} if they have the same grandfather but they are not brothers.

A backward path in the tree is a path going from one site to one of its ancestor and containing all intermediate ancestors.
A forward (finite) path  is a backward path read in the other direction.

\bigskip
The distance between two trees $\GA$ and $\GB$ is $2^{-N(\GA,\GB)}$ where $N(\GA,\GB)$ is the minimal integer $n$ such that $\GA_{\om}\neq \GB_{\om}$ and $|\om|=n$.

In other words, $d(\GA,\GB)=2^{-n}$ means that $\GA$ and $\GB$ have different root if $n=0$, and $\GA_{\om}$ and $\GB_{\om}$ do coincide  for every $\om$, such that $|\om|\le n$ and for at least one $\om$ with $|\om|=n$,  one  of the followers of $\GA_{\om}$  is different to the same follower for $\GB_{\om}$.

Note that the space of trees $\ST$ is compact (for the metric we introduced) as a product of compact spaces. The subset of trees with root equal to 0 (or 1) is also compact as a closed set included into a compact set.


\subsubsection{Canonical dynamics on trees}

There is a natural $\FFF$-action on trees defined as follows:

 For any given $x \in \FF_2^+$  and $\GA \in  \FFF \to \{0, 1\} $ take
$$
      (T_a \GA)_x : = \GA_{ax}    \qquad \text{and} \qquad  (T_b \GA)_x : = \GA_{bx}
$$

For the definition of  distance  one can easily see that we  have:
$$
   d(T_{\om}(\GA), T_{\om}(\GB)) \leq 2 d(\GA, \GB)   \qquad \text{for $\om = a, b$}
$$
and thus the maps $T_a$ and $T_b$ are Lipschitz-continuous.

\begin{definition}
\label{def-patch}
Let $\GA$ be a tree. A patch (of size $n\ge 1$) is any ball $B(T_{\om}\GA,2^{-n})$ where $\om\in \FFF$
\end{definition}

\begin{definition}
\label{def-period}
A tree $\GA$ is called pre-periodic if there exists a finite set of trees $\CO:=\{\GA(1),\GA(2),\ldots, \GA(n)\}$ satisfying
\begin{enumerate}
\item it contains $\GA$.
\item It is $T_{a}$ and $T_{b}$-invariant.
\end{enumerate}
The tree $\GA$ is said to be periodic if $T_{a}(\CO)\cup T_{b}(\CO)=\CO$.
\end{definition}
Equivalently,  $\GA$ is periodic and $\CO:=\{\GA(1),\GA(2),\ldots, \GA(n)\}$ is its orbit, means that the if we consider the oriented graph whose states are the elements of $\CO$ and with edges representing images by $T_{a}$ or $T_{b}$, then, each vertex is the initial point for two edges, and the end point for at least one edge.

\begin{notation}
For simplicity we shall set $T^{-1}(\GA)$ for $T_{a}^{-1}(\GA)\cup T_{b}^{-1}(\GA)$.
\end{notation}

\subsubsection{Substreetution}

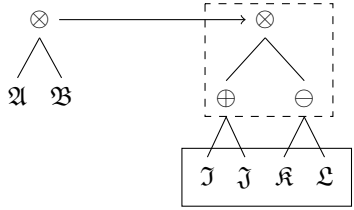
\begin{figure}[h]
\begin{tikzpicture}
\begin{scope}
\Tree [.\node (0) {$\otimes$} ; [.\node {$\GA$} ;  ] [.\node (r) {$\GB$}; ] ]
\end{scope}
\begin{scope}[xshift=3cm]
\Tree [.\node (1) {$\otimes$} ; [.\node(2) {$\oplus$} ; \node(a){$\GI$};  \node(b){$\GJ$};  ] [.\node (r) {$\ominus$}; \node(c){$\GK$}; \node(d){$\GL$}; ] ]
\node[draw, fit=(a) (b) (c) (d) ]{};
\end{scope}
\draw[->] (0)..controls +(east:2) and +(west:1)..(1);
\draw[dashed] (2.2,0.3) rectangle ++(1.7,-1.5);

\end{tikzpicture}
\caption{A market substreetution}
\label{Fig1-substi}
\end{figure}

\begin{definition}
\label{def-sub}
A {\bf substreetution}\footnote{Actually this is a \emph{constant length 2} substitution.} on trees is a map $H$ on the set of configurations defined by concatenation as follows:
\begin{enumerate}
\item $H$ maps each site to a truple (actually a root with two followers), the value depending only on the value of the digit at the site. See Figure \ref{Fig1-substi} with the box with dashline.
\item $H$ connects images of subtrees (followers) as indicated on Figure \ref{Fig1-substi}, with $\GI,\GJ,\GK,\GL\in \{H(\GA),H(\GB)\}$.
\end{enumerate}
 The order word $\GI \GJ \GK \GL$ is called the {\bf grammar} of the substreetution and is explained below.

The substreetution is said to be \emph{marked} if $H(0)=\tritree{$i$}{}{}$ and $H(1)=\tritree{1-i}{}{}$, $i=0,1$.
\end{definition}

\medskip

Identifying the trees and their name, and reemploying notation from above, the word $\GI\GJ\GK\GL$ is a word in $\{H(\GA),H(\GB)\}^{4}$, and can canonically be identified with a word in $\{\GA,\GB\}^{4}$. This defines the \emph{grammar} of the substreetution $H$ and indicates the $H$-image of which sub-tree as to be connected on which slot at generation 2.

\eg  The  grammar $ABBA$ means $\GI=\GL=H(\GA)$ and $\GJ=\GK=H(\GB)$.  In this case it is easy to see that the function
$H$ satisfies the following relations: $H T_a = T_{aa} H = T_{bb} H$ and $H T_b = T_{ab} H = T_{ba} H$.

\bigskip
More generally, we let the reader check that
\begin{equation}
\label{eq-renorm}
\{T_{aa}H, T_{ab}H,T_{ba}H,T_{bb}H\}=\{HT_{a},HT_{b}\},
\end{equation}
which is rewritten under the form $\disp\boxed{T^{2}\circ H=H\circ T}$. We call this relation the \emph{renormalization equation}.

We emphasize
that Equality \eqref{eq-renorm} is a natural extension to the $\FFF$-action case of property $\s^{2}\circ H=H\circ \s$ which holds for the shift considered as a $\NN$-action and for constant length-2 substitutions.

{
\bigskip
It has been pointed to us that our procedure could have some similarities with the one introduced in \cite{Damm}. This procedure deals with theoretical computer science. The vocabulary and probably the questions are quite different from the ones in Dynamical System\footnote{at least as used or studied by authors.}.
This makes the understanding of similarities or differences harder. At least, we see a difference as the Damm procedure seems to study finite trees whereas we are interested in infinite trees. We emphasize that our procedure has been defined from the need that the renormalization equality \eqref{eq-renorm} does hold. }

\subsubsection{Results}

The two first results describe some property of one specific substreetution.

\begin{maintheorem}\label{MainThA}
Let $H$ denote the substreetution given by $\disp 0\mapsto$ \begin{forest}
[$0$[$1$][$0$]]
\end{forest} and $\disp 1\mapsto $\begin{forest}
[$1$[$1$][$0$]]
\end{forest}, equipped with the grammar BBAB.
Then there exists a unique fixed point $\GJ$ with root 0. The closure of its orbit is a minimal dynamical system $X$ and is not periodic.

There also exists a unique fixed point $\GJ'$ with root 1. It coincides with $\GJ$ except at the root.
\end{maintheorem}

{ \noindent
As $\s$ is continuous and $X$ is compact, $\s(X)\subset X$ and minimality yield that $\s(X)=X$. Hence $\s$ (restricted to $X$) is onto, which is {\it a priori} not obvious. This makes sense to study the backward complexity in $X$.
Furthermore, it is neither clear that $\GJ'$ belongs to $X$. This actually holds, see Corollary \ref{coro-jprimeinx} }

\begin{maintheorem}\label{MainThB}
With the same $H$.
Let $\GA$ be in $X$. Set $p(n,\GA):=\#\{\GB\in X,\ T_{\om}(\GB)=\GA, |\om|=n\}$.
Then for every $n$,
$$p(n,\GA)\le 3^{n}.$$
\end{maintheorem}

The bound in Theorem \ref{MainThB} is not sharp at all. What is remarkable it is that it exists.
We emphasize that in the whole space of trees $\CA^{\FFF}$,  for any tree $\GA$, $T^{-1}(\GA)$ has cardinality of the continuum. Furthermore, we conjecture a result similar to what holds for the tree-pressure in \cite{Przytycki-pressure}: for \emph{most} of the elements $\GA$ in $X$,
$$\lim_{\ninf}\frac1n\log p(n,\GA)=\sup_{\GB\in X}\lim_{\ninf}\frac1n\log p({n,\GB}),$$
and both quantities exist.

Re-employing vocabulary from 1-$d$-ergodic theory,
Balls of size $2^{-n}$ form a partition of $X$. We denote by $R_{n}$ its cardinality. $R_{n}$ is the number of patches of size $n$ that we see in $\GJ$. We conjecture that
$$\disp\lim_{n\to+\8}\frac1n\log R_{n}=\sup_{\GB\in X}\lim_{\ninf}\frac1n\log p({n,\GB})$$
holds.

\bigskip
Next results show that substreetutions may not generate a minimal set.
And the last result shows that existence of invariant measures on a periodic set of trees  does not always hold.

\begin{theorem}
\label{th-nonminimal}
The substreetution given by $\disp 0\mapsto$ \begin{forest}
[$0$[$0$][$1$]]
\end{forest} and $\disp 1\mapsto $\begin{forest}
[$1$[$1$][$0$]]
\end{forest}, equipped with the grammar ABBA is non-minimal.
\end{theorem}

\begin{theorem}\label{th-perio-nomeas}
There exist periodic trees with no invariant measure on their orbit.
\end{theorem}

\subsection{More general substreetutions}
Once we agree that the way to define a substitution on trees is to consider that equality in the renormalization equation \eqref{eq-renorm} holds for sets, then it is very easy to define more general substreetutions.

Actually, $H$ can be defined by:
\begin{itemize}
\item Fixing the image of each ``color'' as a finite tree.
\item Fixing the grammar for each \emph{slot} at the end of each $H(\otimes)$.
\end{itemize}
This may be defined for colored trees involving a greater alphabet (more than 2 colors) and even not necessarily binary (or even regular) trees. Of course, the grammar has to consider possible irregularities in the number of followers.

\eg
We can consider a substreetution on binary 2-colored trees defined by $H\left(\tritree{$0$}{$A$}{$B$}\right)=
\begin{forest}
[$0$[$0$[$H(A)$][$H(B)$]][$H(B)$]]
\end{forest}$ and $H\left(\tritree{$1$}{$A$}{$B$}\right)=\heptatree{$1$}{$0$}{$0$}{$H(A)$}{$H(\reflectbox{B})$}{$H(B)$}{H(A)}$.

\bigskip
In this paper, we will mostly consider constant length-2 substreetutions. First of all because these objects are already complicated and there are many situations that we have not yet studied.
For constant-length 2 substitutions on alphabet with 2 letters there are $2^{3}=8$ possibilities. For constant length-2 substreetutions on regular binary trees there are 128 possibilities. Even if these numbers may be diminished due to symmetries, the number of possibilities for trees is much higher than the number of possibilities for the line.

Furthermore, we will deeply use the constant-length structure to define tools (as the \emph{source and the $\chi$-procedure}). These tools could be defined for general substreetution but would be extremely more complicated and harder to work with.

\subsection{Plan of the paper}
In Section \ref{sec-gene} we first give some general results on subtreetution, as existence of fixed points. Then, we prove some points of Theorem \ref{MainThA}: existence of two fixed points that differ only from the root.

\medskip\noindent
In Section \ref{sec-geoJ} we give a precise description of each line of the fixed point $\GJ$. For that we use the notion of source which corresponds to unsubstitution in the classical.

\medskip\noindent
Section \ref{sec-mainA} is devoted to the proof of the remaining points in Main Theorem \ref{MainThA}: that $X$ is minimal and that $\GJ$ is not periodic.

Minimality follows from  \cite{JdeVries}. We show that $\GJ$ is almost periodic. The proof that $\GJ$ is not periodic is done by contradiction. Assuming it is periodic, we use the Perron-Frobenius theorem on positive matrices to show that the proportion of 1's along the lines of $\GJ$ should have finitely accumulation points that are all positive. Contradiction comes from the fact that the proportion of 1's goes to 0 along special subsequences of lines.

\medskip\noindent
Section \ref{sec-rigidity} explores several symmetries and rigidities of the configuration $\GJ$ and its images. This allows to determines types (namely odd and even) of trees in $X$.

We also show that knowing one subtree in $\GA\in X$ allows to  to recover his brother and even his cousins.

\medskip\noindent
Section \ref{sec-mainB} is devoted to the proof of Main Theorem \ref{MainThB}. For any type of tree defined above, we compute the number of possible preimages in $X$ and give a precise description of them.

\medskip\noindent
Finally, in Section \ref{sec-Appendix} we show Theorems  \ref{th-nonminimal} and
\ref{th-perio-nomeas}. We also give a funny example of substreetution generating
the well-known Thue-Morse sequence along the paths.

We give a picture of $\GJ$ and show how to get a colored quasi-periodic tiling of $\mathbb{D}^{2}$ using elements of $X$.


\section{Some more general facts on (constant length-2) substreetutions}\label{sec-gene}


\subsection{Source and Unsubstreetuted tree}\label{ss-unsub}
We consider a \emph{constant length 2} substreetution as defined in Def. \ref{def-sub}.

We remind the renormalization equation $T^{2}\circ H=H\circ T$. From this we get:

\begin{lemma}
\label{lem-source}
Let $H$ be a substreetution.
For every site $\omega$ with even length, there exists a word $s(\om)$ of length $|\om|/2$ such that
$$T_{\om}\circ H=H\circ T_{s(\om)}.$$
The site $s(\om)$ is called the \emph{source} of the site $\om$.
\end{lemma}
\begin{proof}
Direct application of renormalization equation by induction.
\end{proof}

We emphasize that a kind of converse is true. Given any word $\om$, there exists $\om'$ with $|\om'|=2|\om|$ such that $T_{\om'}\circ H=H\circ T_{\om}$. Nevertheless, such an $\om'$ is not necessarily unique. It is entirely determined by the grammar.

\subsection{Fixed points for marked substitution}

We start with a basic observation.
\begin{lemma}
\label{lem-subsinj}
A marked substreetution is one-to-one.
\end{lemma}
\begin{proof}
Let $\GA$ and $\GB$ be two trees coinciding up to generation $n-1$. Let $\om$, with $|\om|=n$ such that $\GA_{\om}\neq\GB_{\om}$. Let $\om'$ be such that $s(\om')=\om$. By definition, $T_{\om}(\GA)$ and $T_{\om}(\GB)$ have different roots. Because $H$ is marked, $H\circ T_{\om}(\GA)=T_{\om'}(H(\GA))$ and $H\circ T_{\om}(\GB)=T_{\om'}(H(\GB))$ have different roots. Therefore, $H(\GA)\neq H(\GB)$.
\end{proof}

Assume that $\GA$ is a fixed-point for a substreetution $H$. Let $\om$ be a site with even length. Set $\om'=s(\om)$. Then, $T_{\om}(\GA)=H(T_{\om'}(\GA))$. We say that $T_{\om'}(\GA)$ is the \emph{unsubstreetuted tree} from $T_{\om}(\GA)$. By lemmas \ref{lem-source} and \ref{lem-subsinj} it is well defined and uniquely determined. However, because the source operation is not one-to-one, it may be that several sites have the same unsubstreetuted tree. This hold if and only if the subtrees with roots these sites are equal.

Now, we show  that  a well-known result on classical substitutions holds for substreetution.

\begin{proposition}
\label{prop-pointfixeH}
If $H(0)=\arb{0}$ (resp. $H(1)=\arb{1}$) then, there is a unique fixed point for $H$ starting with $0$ (resp. 1). They will respectively be denoted by $\afz$ and $\afu$.
\end{proposition}
\begin{proof}
Let us assume that  $H(0)=(0,.,.)$ holds,  the other case being totally symmetric.

By definition of $H$, if $\GA$ and $\GB$ are two trees coinciding up to generation $n$, with $n\ge 1$, then $H(\GA)$ and $H(\GB)$ do coincide up to generation $2n$ because each site produces 2 generations of new sites. Therefore,
$$d(H(\GA),H(\GB))\le d^{2}(\GA,\GB)\le \frac12d(\GA,\GB).$$
Hence, $H$ acts as a contraction on the compact set of trees whose root is 0. It admits a unique fixed point on this set.
\end{proof}

Because the set of trees is compact, thus is an Hausdorff space, the fixed point $\afz$ is obtained as the limit of $H^{n}(\GA)$ as $n$ goes to $+\8$ and where $\GA$ is any tree with root $0$. Furthermore, $\afz$ can also be obtained as $\lim_{\ninf}H^{n}(0)$.

In Theorem \ref{MainThA} we consider the substitution
$\disp 0\mapsto$ \begin{forest}
[$0$[$1$][$0$]]
\end{forest} and $\disp 1\mapsto $\begin{forest}
[$1$[$1$][$0$]]
\end{forest}, equipped with the grammar BBAB.

It fulfills conditions to ensure that there is a fixed point with root $0$ and a fixed point with root $1$. We denote by $\GJ$ the fixed tree with root $0$ and by $\GJ'$ the fixed point with root 1. It is easy to check (by induction) that each odd  rank is a string of $10$'s.

We give here a first description  of these two fixed points:
\begin{proposition}
\label{prop-Jprime}
The two trees $\GJ$ and $\GJ'$ only differ from  the root.
\end{proposition}
\begin{proof}
We recall that $\GJ$ and $\GJ'$ are obtained as $\lim_{\ninf}H^{n}(0)$ and $\lim_{\ninf} H^{n}(1)$.
Hence, the proof is done by induction on the number of iterations for $H$. we show that $H^{n}(0)$ and $H^{n}(1)$ only differ in the root.

Because $H(0)=\tritree{0}{1}{0}$ and $H(1)=\tritree{1}{1}{0}$, then level 1 for $\GJ$   is equal to level 1 for $\GJ'$.

For $n\ge 2$, $H^{n}(0)=H^{n-1}\left(\tritree{0}{1}{0}\right)$. This  finite tree is the tree $\GA^{n-1}:=H^{n-1}(0)$ where we glue at the extremities the trees $H^{n-1}(0)$ or $H^{n-1}(1)$. The gluing only depends on the grammar rules.

Similarly, the finite tree $H^{n}(1)=H^{n-1}\left(\tritree{1}{1}{0}\right)$ is equal to the tree $\GB^{n-1}:=H^{n-1}(0)$ where we glue at the extremities the trees $H^{n-1}(0)$ or $H^{n-1}(1)$. The gluing for that tree respects the same rules that for $H^{n}(0)$.

By assumption, $\GA$ and $\GB$ only differ in the root.
This finishes the proof by induction.
\end{proof}

\bigskip
In the following, we write $\CO(\GJ)$ for $\disp \{T_{\om}(\GJ),\ \om\in\FFF\}$.
 We remind that $X=\ol {\CO(\GJ)}$.

\section{Description of the Jacaranda tree $\GJ$}\label{sec-geoJ}
In all the following, except if it is especially mentioned, we shall only consider the substitution from Theorem \ref{MainThA}:
$\disp 0\mapsto$ \begin{forest}
[$0$[$1$][$0$]]
\end{forest} and $\disp 1\mapsto $\begin{forest}
[$1$[$1$][$0$]]
\end{forest}, equipped with the grammar BBAB. In this case the renormalization equations take the form
\begin{equation}\label{eq-renormBBAB}
 T_a T_a H = T_b T_a H = T_b T_b H = T_b, \;\; T_a T_b H = H T_a
\end{equation}
$\GJ$ and $\GJ'$ are the two fixed points (for $H$) respectively with root 0 and 1. $\GJ$ is also referred to as the \emph{Jacaranda tree}.

\bigskip
As it is said above, each odd line in $\GJ$ is a concatenation of 10's. The goal of this section is first to describe each even line (see Prop. \ref{prop-descrifixed1}) and to get an estimation for the number of 1's on even lines. This later result will be used to show that  $\GJ$ is not periodic (nor pre-periodic).


\subsection{Inverse of the source function for $\GJ$}

  As defined previously in the lemma \ref{lem-source}, the source function $s$ is defined by $$
  T_{\omega} H = H T_{s(\omega)}.
$$

 Here we want to describe more precisely the  function $s$; the map $H$ has the property that for any
 given configuration $\GC$ we get $H(\GC)_e = \GC_e$. We also have that
 $$
   aaH = baH = bbH = Hb  \;\; \text{and} \;\; abH = H a
 $$
 (where, for example $abH$,  is just a short notation for $T_a T_b H = T_{ab} H $).

If we consider $\GJ$,  the fixed point for $H$ whose root is zero, then we have
$$
    (a_1 a_2 \ldots a_k \GJ)_e = \GJ_{a_k \ldots a_2 a_1}
$$

 We also get, using equation \eqref{eq-renormBBAB},
$$
  \GJ_{ba} = (ab\GJ)_e = (ab H(\GJ))_e = (Ha(\GJ))_e = (a\GJ)_e = \GJ_a = \GJ_{s(ba)}
$$
Hence the source of $ba$ is $a$. In a similar way we can show that the sources of $aa, ab, bb$ are $b$:
$$
\begin{cases}
s(aa)=s(ab)=s(bb)=b\\
s(ba)=a,\\
\end{cases}
$$

More generally we get:
\begin{lemma}\label{lem-propsource} Given two words $P$ and $Q$ with even length then
  $s(P Q ) = s(P) s(Q)$
\end{lemma}
\begin{proof}
The proof uses the renormalization equation \eqref{eq-renormBBAB}: we have
$$
 \GJ_{p_1 p_2 p_3 p_4 \ldots p_{2n-1}p_{2n}} = (p_{2n} \ldots p_4 p_3 p_2 p_1 \GJ)_e = (p_{2n} \ldots p_4 p_3 p_2 p_1 H(\GJ))_e =
$$
$$
= (p_{2n} \ldots p_4 p_3 H s(p_1 p_2) \GJ)_e = (p_{2n} p_{2n-1} \ldots H s(p_3 p_4) s(p_1 p_2) \GJ)_e =
$$
$$
 \ldots = H(s(p_{2n-1} p_{2n}) \ldots s(p_3 p_4) s(p_1 p_2) \GJ  )_e = \GJ_{s(p_1 p_2) s(p_3 p_4) \ldots s(p_{2n-1} p_{2n}) }
$$
and so we obtain
$$
  s(p_1 p_2 p_3 p_4 \ldots p_{2n} ) = s(p_1 p_2) s(p_3 p_4) \ldots s(p_{2n-1} p_{2n})
$$
Hence,  if $P$ and $Q$ are two words with even length  we then have
$s(P Q ) = s(P) s(Q)$, as claimed.
\end{proof}

Now we want to see the inverse image of $s$, that we denote by $\theta$, and it is seen as a multivalued function
that maps words to  sets of words.
It is defined in such way that $ s( \theta(p)) = \{p\}$  for any word $p$ and $\theta (s (P)) \ni P$ for any
even length word $P$; in the definition of $\theta$ we also set $\theta(\emptyset) = \emptyset $ and $\theta(e) = \{e\}$.

The main property of the function $\theta$ is the following:
\begin{lemma}\label{lem-proptheta} Given the words $p$ and  $q$ then
  $\theta(pq) = \theta(p) \theta(q)$
  (where the product above is the set of words obtained by concatenating any word in $\theta(p)$ with any another word in $\theta(q)$).
\end{lemma}
\begin{proof}
($\theta(p)\theta(q) \subseteq \theta(p q)$):
 Let $R \in \theta(p) \theta(q)$; then $R = P Q$ where $P \in \theta(p)$ and $Q \in \theta(q)$.
Hence $s(P) = p$ and $s(Q) = q$ and Lemma \ref{lem-propsource} yields $s(R) = s(PQ) = s(P) s(Q) = p q $, showing that $pq$ is the source for
$R$. Therefore $R \in \theta(pq)$ holds.

($\theta(pq) \subseteq \theta(p) \theta(q)$):
Let $R \in \theta(pq)$; hence $s(R) = pq$. Since $R$ is a word with even length (because it is in the image of $\theta$)
  we can write $R = A B $ with $|A| = 2 |p|$ and $|B| = 2 |q|$. Then $s(R) = s(A B) = s(A) s(B)$. This shows that $s(A)$ is the prefix of $s(R)=pq$ with length $|p|$. Similarly $s(B)$ is the suffix of $pq$ with length $|q|$. Hence,
  $s(A) = p$ and  $s(B)= q$. This implies that $A \in \theta(p)$ and $B \in \theta(q)$, showing that $R = A B \in \theta(p) \theta(q)$
  as claimed, concluding the proof.
\end{proof}

A simple recurrence yields
\begin{equation}
\label{eq-propthetegene}
\theta(p_1 p_2 \ldots p_k) = \theta(p_1) \theta(p_2) \ldots \theta(p_k).
\end{equation}
It is also easy to see from the definition that $\theta(a) = ba$ and $\theta(b) = \{aa, ab, bb  \}$; hence, using the property
above we can obtain $\theta$ of any given word.

\subsection{The $\chi$-procedure for even lines}
Odd lines for $\GJ$ and $\GJ'$ are just string of $10$'s. We want to get a description for even lines. We start with the last even line for finite trees.

\subsubsection{Word at last even line for image of finite colored trees }
Here, we want to describe how the word at last line for a colored tree generates the word at last but one line for its image by $H$. This procedure is recursive.

\begin{definition}
\label{def-familyblock}
A family-block in a colored tree  is the collection of all $n$-descendants of a site $\om$.
\end{definition}
Note that 2 brothers and 4 cousins form a family block. A family block has cardinality $2^{n}$, where $n$ is such that the first common ancestor is their $n$-ancestor.

\begin{definition}
\label{def-bottomword}
Let $\GA$ be a finite colored tree whose last line is a full family block. Then, the word at the last line is called the bottom-word of $\GA$. It is denoted by $\bw(\GA)$. The word at line last but one is called the vice-bottom-word. It is denoted by $\vbw(\GA)$.
\end{definition}
Next lemma is a direct consequence of the definition of $H$:
\begin{lemma}
\label{lem-bwh}
Let $\GA$ be a finite colored tree whose last line is a full family block. Set $\GA=\tritree{$\otimes$}{$\GL$}{$\GR$}$. Then
$$\vbw(H(\GA))=\vbw(H(\GR))\vbw(H(\GR))\vbw(H(\GL))\vbw(H(\GR)).$$
\end{lemma}
Note that equality $\bw(H(\GA))=\bw(H(\GR))\bw(H(\GR))\bw(H(\GL))\bw(H(\GR))$ also holds but is not interesting as the bottom-word is only composed by a chain of 10's.

{  Now, in order to introduce $\chi$, we notice that
 we can represent any configuration $\GC \in \{0, 1\}^{\FFF} $ by means of the sets of addresses
 of the 1's in each of its lines (it is clear that this can also be done
 using zeros, our choice of ``1'' here being completely arbitrary).   Let us define
 the set $W_n  = \{a, b\}^n$ of words of length $n \geq 1$ over the alphabet $\{a, b\}$, $W_0$ being the singleton $\{e\}$
 that represents the root.

 Given a configuration $\GC$ we have, at the line $l \geq 0$, a word of length $2^l$ on the alphabet $\{0, 1\}$ whose
 addresses are words on the alphabet $\{a, b\}$, with length $l$; hence the addresses corresponding to the
 line $l$ are in $W_l$. Then we define $1_l(\GC)$ as the set of addresses $\omega \in W_l$ where
 $\GC_{\omega} = 1$.

For example $\mathbf{1}_2(0010) = \{ ba \}$ and $\mathbf{1}_3(01000001) = \{ aab, bbb \}$;
in the case
of words with length one (corresponding to the root of the configuration) we  write
$\mathbf{1}_0(0) = \emptyset$ and $\mathbf{1}_0(1) = e$.

Since the alphabet is $\{0, 1\}$ we have a one to one correspondence between a word $P \in \{0, 1\}^{2^l}$ (with length $2^{l}$)
and the set $\mathbf{1}_l(P) \subset W_l$.

\medskip

In particular, this allows us to define a function $\chi$ that maps words of $\{0, 1\}^{2^l}$ to words of $\{0, 1 \}^{2^{2l}}$:
given a word $P  \in \{0, 1\}^{2^l}$ we take the corresponding addresses of the symbols ``1'' (say, $\mathbf{1}_l(P)$),
take their image under $\theta$, that gives a set of addresses in $\{a, b\}$ with length $2l$ and consider then as the
addresses of the symbols ``1'' of a new word $\chi(P)$ of length $2^{2l}$. More formally:
\begin{definition}
$\chi$ is the function  such that
$$
   \mathbf{1}_{2l}(\chi(P)) = \theta(\mathbf{1}_l(P))
$$
where $P$ is a word of length $2^l$ over the alphabet $\{0, 1\}$.
\end{definition}
\begin{remark}
 In the expression above, as usual, $\theta(A) = \{  \theta(a) \; \forall a \in A \}$.
\end{remark}

From this definition we get that the function $\chi$ can be described recursively:
\begin{lemma}
\label{lem-chi}
We have that $ \chi(0) =0, \chi(1) = 1$
and, for any given word $W_1 W_2$ with length $2^l$, we have
$$
   \chi(W_1 W_2) = \chi(W_2) \chi(W_2) \chi(W_1) \chi(W_2)
$$
\end{lemma}
\begin{proof}
We have
$$
 \mathbf{1}_0(\chi(0)) = \theta(\mathbf{1}_0(0)) = \theta(\emptyset) = \emptyset
$$
Hence, $\chi(0) = 0$.
Similarly,
$
 \mathbf{1}_0(\chi(1)) = \theta(\mathbf{1}_0(1)) = \theta(e) = e
$
and so $\chi(1) = 1$.

Now
$$
\mathbf{1}_{2l}(\chi(W_1 W_2)) = \theta(\mathbf{1}_l(W_1 W_2))
$$
But $\mathbf{1}_l(W_1 W_2) = \{a \mathbf{1}_{l-1}(W_1), b \mathbf{1}_{l-1}(W_2)   \}$ and then
$$
 \theta(\mathbf{1}_l(W_1 W_2)) = \theta( \{a \mathbf{1}_{l-1}(W_1), b \mathbf{1}_{l-1}(W_2)   \} ) =
$$
$$
 = \{  ba \theta(\mathbf{1}_{l-1}(W_1)), aa \theta(\mathbf{1}_{l-1}(W_2)),
 ab \theta(\mathbf{1}_{l-1}(W_2)),  bb \theta(\mathbf{1}_{l-1}(W_2)) \}   =
$$
$$
 = \{  aa \mathbf{1}_{2l-2}(\chi(W_2)), ab \mathbf{1}_{2l-2}(\chi(W_2)),
 ba \mathbf{1}_{2l-2}(\chi(W_1)), bb \mathbf{1}_{2l-2}(\chi(W_2)) \}  =
$$
$$
   \mathbf{1}_{2l}(\chi(W_1 W_2))
$$
Hence $\chi(W_1 W_2) = \chi(W_2) \chi(W_2)  \chi(W_1) \chi(W_2) $
\end{proof}
}

\begin{tabular}{|p{10cm}|}
\hline
 We emphasize (and let the reader check):
 \begin{enumerate}
\item $\chi(10)=0010$,
\item any $\chi^{u}(10)$ with $u\ge 2$ is a concatenation of blocks 0010 and 0000,
\item $\chi^{u}(10)$ is a word of length $2^{2^{u}}$ in 0's and 1's.
\end{enumerate}\\
\hline
\end{tabular}

\begin{proposition}
\label{prop-vbw-chi}
For any finite colored tree with full last line, say $\GA$, then $\vbw(H(\GA))=\chi(\bw(\GA))$ holds.
\end{proposition}
\begin{proof}
The proof is done by induction on the number of lines of the finite colored tree $\GA$ (with full last line).
If $\GA=0,1$ (that is only line 0), then $H(\GA)=\tritree{$\otimes$}{1}{0}$ with $\otimes=\GA$. Moreover, $\bw(\GA)=\otimes$.
Then,
$$\vbw(H(\GA))=\otimes=\chi(\otimes)=\chi(\bw(\GA))$$

Assume that the result holds for any finite colored tree with full last line with $l$ lines. Let $\GA$ be a finite colored tree with full last line and with $l+1$ lines. Set $\GA=\tritree{$\otimes$}{$\GL$}{$\GR$}$ and  $\bw(\GA)=w$. Then $w=\bw(\GL)\bw(\GR)$  holds with $\GL$ and $\GR$ two colored trees with full last line and with $l$-lines.

Lemma \ref{lem-bwh} and Induction Assumption yield
\begin{eqnarray*}
\vbw(H(\GA))&=& \vbw(H(\GR))\vbw(H(\GR))\vbw(H(\GL))\vbw(H(\GR))\\
&=&  \chi(\bw(\GR))\chi(\bw(\GR))\chi(\bw(\GL))\chi(\bw(\GR))\\
&=& \chi(\bw(\GL)\bw(\GR))=\chi(\bw(\GA)).
\end{eqnarray*}

\end{proof}

{
\subsubsection{Some technical lemmas on $\chi$}

We give some extra properties for $\chi$ that will be used later.
\begin{lemma}
\label{lem-techchi}
If $\om$ is a word with length $2^{l}$ and $\om'$ is the concatenation of $2^{k}$ words $\om$, then $\chi(\om')$ is a concatenation of $\chi(\om)$.
\end{lemma}
\begin{proof}
The proof is done by induction on $k$. The result is obvious if $k=1$. Let us assume it holds for $k$ and let us prove it for $k+1$.

Set $\om'$ the concatenation of $2^{k}$ words $\om$ and $\om''=\om'\om'$.
Then,
$$\chi(\om'')=\chi(\mathbf{\om'}\om')=\chi(\om')\chi(\om')\chi(\mathbf{\om'})\chi(\om').$$
Hence, the property holds for $k+1$ just by using Induction Hypothesis.
\end{proof}

\begin{lemma}
\label{lem-chiades1}
For any $u\ge 1$, the word $\chi^{u}(10)$ has at least one 1 in the second half. The same holds for the first half if $u\ge 2$
\end{lemma}
\begin{proof}
The proof is done by induction.

First, $\chi(10)=0010$. The second half of the word is 10 which contains 1.

Let us assume that the result holds for $\chi^{u}(10)$ and let us prove it for $\chi^{u+1}(10)$. Set $\chi^{u}(10)=W_{1}W_{2}$ with $|W_{1}|=|W_{2}|$.

By definition, $\chi^{u+1}(10)=\chi(\chi^{u}(10))=\chi(W_{1}W_{2})=\chi(W_{2})\chi(W_{2})\chi(W_{1})\chi(W_{2})$. The second half of $\chi^{u+1}(10)$ is the word $\chi(W_{1})\chi(W_{2})$ as the length of the image of a word by $\chi$ does only depend on the length of the word.

By Induction hypothesis, $W_{2}$ contains some 1. By definition of $\chi$, as $\chi(1)=1$, $\chi(W_{2})$ also contains some 1.  This shows that both halves contain some 1.
\end{proof}
}

\subsubsection{Description of words at even line in $\GJ$}

\begin{proposition}
\label{prop-descrifixed1}
For every integer of the form $m=2^{u}(2n+1)$ with $u\ge 0$, the word at generation $m$ in $\GJ$ is a concatenation of blocks  $\chi^{u}(10)$.
\end{proposition}
\begin{proof}
The proof is done by induction on $u$.

For $u=0$, we are looking at an odd line, which is only composed by blocks 10 by definition of $H$.

We assume that this the property holds for any line of the form $2^{u}(2n+1)$ and we prove it for a line of the form $2^{u+1}(2n+1)$.

At line $2^{u}(2n+1)$ we consider all the consecutive and juxtaposed blocks of length $2^{2^{u}}$.
They form a sequence of family blocks with common $2^{u}$-ancestors. More precisely, the line $2^{u}(2n+1)$ is the juxtaposition of the bottom-words of all the finite subtrees in $\GJ$ with root at line $2^{u+1}n$ and with size $2^{u}$.

The tree $\GJ$ is invariant by $H$, hence the line $2^{u+1}(2n+1)$ is obtained from the line $2^{u}(2n+1)$. From Proposition \ref{prop-vbw-chi} and induction hypothesis, Lemma \ref{lem-techchi} shows that the line $2^{u+1}(2n+1)$ is the concatenation of blocks $\chi(\chi^{u}(10))=\chi^{u+1}(10)$.

 \end{proof}

\begin{tabular}{|p{10cm}|}
\hline
 We emphasize:
 \begin{enumerate}
\item odd lines are only composed by concatenations of 10's.
\item Lines in $4\N+2$ are only composed by 0010.
\item Lines in $4\N$ are composed by blocks 0010 and blocks 0000.
\end{enumerate}\\
\hline
\end{tabular}

\subsection{Proportion of 1's on even lines in $\GJ$}

{
Now we want to estimate how many 1's we have in even lines of the Jacaranda tree $\GJ$.

First we start with special even lines.
\begin{proposition}
 The number of 1's on the line $2^n$ is given by $2^{2^n} (1+f^n(1))^{-1}$,
 where $f(x) = x+ x^{-1} + 1$ and $n \geq 0$.
\end{proposition}
\begin{proof}
First, we recall the map $\theta$, defined  such that $\theta(a) = ba, \theta(b) = \{ aa,  ab, bb   \} $ and
  $\theta(p q) = \theta(p) \theta(q)$ (meaning the concatenation
of any word of the set $\theta(p)$ with any word of the set $\theta(q)$).

Now consider the line $1= 2^0$; since it is an odd line, the element that corresponds to
$1$ is on the site $a$. We define $P_0$ as the number of sites with $1$ and $Q_0$ as the number
of sites with $0$ on the line $2^0$; hence $P_0 = Q_0 = 1$.
This line is the source of line $2 = 2^1$ and so, to find the sites that correspond
to $1$, we need to find all
the sites whose source is $a$. But they are given exactly by $\theta(a) = ba$.

Then on the line $2^n$ we have that the sites with $1$ are given by $\theta^n(a)$ (and the sites
with $0$ are given by $\theta(b)$). Denoting by $P_n$ (resp. $Q_n$) the cardinality of $\theta^n(a)$ (resp. $\theta^n(b)$)
we have that $P_n + Q_n = 2^{2^n}$, the length of the line $2^n$.

To obtain the line $2^{n+1}$, whose source is $2^n$, we use the fact that the ones are on the positions
$\theta(\theta^n(a)) = \theta^n(\theta(a)) = \theta^n(ba) = \theta^n(b) \theta^n(a)$; hence the sites
of the ones are the concatenation
of any site of $\theta^n(b)$ with any from $\theta^n(a)$.
This implies that the number of $1$'s
is given recursively by $P_{n+1} = P_n Q_n$.   To obtain $Q_{n+1}$ we just remind that this number
is the number of elements of the line $2^{n+1}$ minus the number of elements that are equal to $1$:
$Q_{n+1} = (P_n + Q_n)^2 - P_{n+1} = (P_n + Q_n)^2 - P_n Q_n$.

For any $n \geq 1$ we have then the recursion
$$
  P_{n+1} = P_{n} Q_{n}   \qquad \text{and} \qquad Q_{n+1} = (P_{n} + Q_n)^2 - P_{n} Q_n
$$

Writing $R_n = Q_n/P_n$ we obtain that
$$
  R_{n+1} =  \frac{Q_{n+1}}{P_{n+1}} = \frac{(Q_n+P_n)^2- P_n Q_n}{P_n Q_n}   =
$$
$$
 = \frac{Q_n^2 + P_n^2 + P_n Q_n}{ P_n Q_n} = \frac{Q_n}{P_n}+ \frac{P_n}{Q_n} + 1 = R_n + \frac{1}{R_n} +1    = f(R_n)
$$
where $f(x) = x+ \frac{1}{x} +1$ and $R_1 = 3$.

Hence
$$
   2^{2^n} = P_n + Q_n = P_n + R_n P_n \Rightarrow P_n  = \frac{2^{2^n}}{1+f^{n-1}(3)} =
    \frac{2^{2^n}}{1+f^{n}(1)}
$$
as claimed.
\end{proof}
\begin{corollary}\label{coro-cocorico}
 Using the same notation of the previous proposition we have the following:
\begin{enumerate}
\item  The proportion of 1's in the line $2^n$ is given by $P_n/( 2^{2^n}) =  (1+f^n(1))^{-1} $,
that is a strictly decreasing function of $n$ and  this proportion goes to zero as $n$ goes to infinity.
\item  The proportion of 1's in the lines $m=2^u(2n+1)$ is given by $ (1+f^u(1))^{-1} $.
\end{enumerate}
\end{corollary}
\begin{proof}
{\em (1):} notice that $f(x) = x+\frac{1}{x} + 1 > x+1$; hence
 we have that $f^k(1) \geq 1 + k$, showing that $f^k(1) \to \infty$ as $k$ goes to infinity and so
$$
    \frac{P_n}{2^{2^n}} = \frac{1}{1+f^{n}(1)}  \to 0 \;\; \text{as $n \to \infty$}
$$
 It is also clear that $f^n(1) = f(f^{n-1}(1)) = f^{n-1}(1) + \frac{1}{f^{n-1}(1)} + 1 > f^{n-1}(1)$ showing
 that $f^n(1)$ increases with $n$ and so the proportion of 1's decreases as $n$ increases.

{\em  (2):} as a consequence of proposition
 \ref{prop-descrifixed1} we have that the lines $m= 2^u(2n+1)$ are a concatenation of blocs
 $\chi^u(10)$; since in each one of those blocs the proportion of 1's is given by $ (1+f^u(1))^{-1} $
 then the proportion of 1's at the line $m$ is also given by   $ (1+f^u(1))^{-1} $.
\end{proof}

A very useful consequence of the results above is the following:
\begin{corollary}\label{coro-cocorico2}
 A concatenation of blocs $\chi^u(10)$ is the same as a concatenation of blocs
 $\chi^v(10)$ if and only if $u = v$.
\end{corollary}
\begin{proof}
 One implication is automatic, let us concentrate on the other. Suppose $v > u$ (the case
 $v < u$ being completely similar); now take a concatenation of blocs $\chi^{u}(10)$. The proportion
 of 1's in this word is the same as in $\chi^u(10)$, that is larger than the proportion of 1's in
 $\chi^v(10)$. But in the concatenation of $\chi^v(10)$ the proportion of 1's is the same as
 in $\chi^v(10)$, showing that both blocs are distinct when $v \neq u$.
\end{proof}


\section{Proof that $X$ is minimal}\label{sec-mainA}

\subsection{The fixed point is minimal}

The Grammar BBAB yields $T_aT_aH=T_bT_aH=T_b T_bH=H T_b$ and $T_aT_bH=H T_a$. Thus the source is given by
$$\begin{cases}
s(aa)=s(ab)=s(bb)=b\\
s(ba)=a,\\
\end{cases}$$
and more generally, $s(p_{1}p_{2}\ldots p_{2n})=s(p_{1}p_{2})\ldots s(p_{2n-1}p_{2n})$.

\begin{proposition}
\label{prop-sourcelimite1}
There exists $N$ such that for any position $\omega=\om_{1}\ldots \om_{n}$, for any word $\om_{n+1}\ldots \om_{n+N}$, there exists $k\le N$ such that $n+k$ is even and $\GJ_{\om_{1}\ldots \om_{n+k}}=0$.
\end{proposition}
\begin{proof}
The result is obvious if $n$ is even and $\GJ_{\om_{1}\ldots \om_{n}}=0$. For the rest of the proof, we assume that $\GJ_{\om_{1}\ldots \om_{n}}=1$ holds.

Every odd line is followed by an even one.  Hence, we can only prove the proposition for $n$ even, the result for odd line is an immediate consequence of the result for even $n$.

For the rest of proof we consider that $n$ is even and that the digit at positon ${\om_{1}\ldots \om_{n}}$ is $1$.

\begin{enumerate}
\item The case $n\in 8\N$. Let us set $n=8n_{1}$. Because $H$ is marked, the source of the digit at position ${\om_{1}\ldots \om_{n}}$ is $1$ and is an even line. Its source is at line $2n_{1}$ and is again $1$, and the source of this last digit is at line $n_{1}$ and is $1$.
Hence, we have
 line $n_{1}$, $1$ and

\begin{figure}[h]
\begin{tikzpicture}
\begin{scope}
\Tree [.\node (level-1) {$1$} ; $1$ $0$
]
\draw (level-1) node[xshift=1cm] {level $2n_{1}$};
\end{scope}

\begin{scope}[xshift=3cm]
\Tree [.\node (level0) {$1$} ; [.\node(a){$1$} ; \node(b){$0$};  \node(c){$0$};  ] [.\node(d) [draw]{$0$}; $1$ $0$ ] ]
\node[draw, fit=(a) (b) (c)]{};
\draw (level0) node[xshift=1cm] {level $4n_{1}$} ;
\end{scope}

\begin{scope}[xshift=6cm]
\Tree [.\node (level+1) {$1$} ; [.\node {$1$} ; \node(e) [draw] {$0$};  \node [draw] {$0$};  ] [.\node {$0$}; [.\node (0p) {$1$} ; [.\node (1p) {$1$}; \node(3p){$0$}; \node(4p){$0$};
]
[.\node (2p) {$0$} ; \node(5p){$0$}; \node(6p){$0$};
]
]
  \node [draw] {$0$}; ] ]
\draw (level+1) node[xshift=1cm] {level $8n_{1}$};
\node[draw, fit=(0p) (1p) (2p) (3p) (4p) (5p) (6p)]{};
\end{scope}
\draw[semithick,->] (d)..controls +(east:1) and +(west:1)..(e);
\draw[semithick,->] (a)..controls +(west:5) and +(south west:5)..(0p);
\end{tikzpicture}
\caption{Different parts of the tree $\GJ$}
\label{fig-arbre1}
\end{figure}
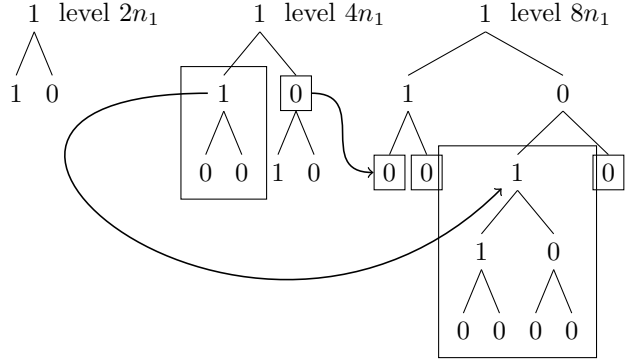

On Figure \ref{fig-arbre1}, site $\boxed{0}$ at level $4n_{1}+1$ is mapped by $H$ on sites $\boxed{0}$ at level $8n_{1}+2$. The  subtree in frame at level $4n_{1}+1$ is mapped on the subtree at level $8n_{1}+2$ in frame (+ children).
At level $8n_{1}+2$ three grand-children are $0$. At level $8n_{1}+4$, all the grand-children from the site with digit $1$ at level $8n_{1}+2$ are $0$.

Consequently, for any $\al$
 with $|\al|\ge 4$, there exists a digit in $\GJ$ at position $\om_{1}\ldots\om_{n}\al_{1}\ldots \al_{k}$ with $k\le 4$ whose value is $0$ and such that $n+k$ is even.

 \item Other cases. Let us set $n=8n_{1}+j$with $j=2,4,6$. Set $\al$ be any word with length at least 10. Then, $\om_{1}\ldots\om_{n}\al_{1}\ldots\al_{8-j}$ has length in $8\N$. Either there is some $0$ at even position from position $\om_{1}\ldots\om_{n}$ to position $\om_{1}\ldots\om_{n}\al_{1}\ldots\al_{8-j}$, or we are back to the previous case.

\end{enumerate}

\end{proof}

{ In the following for $\om\in\FFF=\om_{0}\ldots \om_{n}$, $\wt\om$ stands for $\om_{n}\ldots \om_{0}$.  }

\begin{proposition}
\label{prop-sourcelimite2}
There exists $N_{0}$ such that for any  $\omega=\om_{1}\ldots \om_{n}\in \FFF$,
there exists $k\le N_{0}$  and $\al\in \FFF$ with $|\al|\le k$ such that $T_{\wt\om}(\GJ)$ belongs to $T_{\al}(B(\GJ,1))$.
\end{proposition}
\begin{proof}
From Proposition \ref{prop-sourcelimite1}, path from position $\om':=\om_{1}\ldots\om_{n-N}$ to $\om_{1}\ldots\om_{n}$ contains some $0$ at even level. Having digit equal to $0$ at even level implies being in $B(\GJ,1)$. Thus, result holds with $N_{0}$ equal to $N$ from Proposition \ref{prop-sourcelimite1}.
\end{proof}

\begin{proposition}
\label{prop-minimal}
For every $m$, there exists $N_{m}$ such that for any  $\omega=\om_{1}\ldots \om_{n}\in \FFF$,
there exists $k\le N_{m}$  and $\al\in \FFF$ with $|\al|\le k$ such that $T_{\wt\om}(\GJ)$ belongs to $T_{\al}(B(\GJ,2^{-m}))$.
\end{proposition}
\begin{proof}
The proof is done by induction on $m$. For $m=0$ this follows from Proposition \ref{prop-sourcelimite2}.

Let us assume that the proposition holds for $m$ and let us prove it for $m+1$.

Let $\omega=\om_{1}\ldots \om_{n}\in \FFF$. For simplicity we assume it is bigger than $2N_{m}+N_{0}$. Again, we assume that $n$ is even. From Proposition \ref{prop-sourcelimite2}, we can find $k\le N_{0}$ such that $n-k$ is even and
$$\GJ_{\om_{1}\ldots \om_{n-k}}=0.$$
Let us consider the source for this digit. It lives at level $n'=\disp{n-k\over 2}$. Note that $n'$ is bigger or equal to $N_{m}$.

Let $\be=\be_{1},\ldots \be_{n'}\in\FFF$ be its position. Then, by induction hypothesis, there exists $j\le N_{m}$ such that $T_{\be_{n'-j}\ldots \be_{1}}(\GJ)$ belongs to $B(\GJ,2^{-m})$.

By definition $\be_{1}\ldots \be_{n'}=s(\om_{1}\om_2)\ldots s(\om_{n-k-1}\om_{n-k})$ and then
$$\be_{1}\ldots \be_{n'-j}=s(\om_{1}\om_2)\ldots s(\om_{n-k-2j-1}\om_{n-k-2j}).$$

Then we the renormalization equality and apply $T_{\be_{n'-j}\ldots \be_{1}}$.
This yields
\begin{eqnarray*}
H\circ T_{\be_{n'-j}\ldots \be_{1}}&=& H\circ T_{s(\om_{n'-j}\om_{n'-j-1})\ldots s(\om_{2}\om_{1})}\\
&=& T_{\om_{2n'-2j}\om_{2n'-2j-1}\ldots \om_{2}\om_{1}} H(\GJ)\\
&=&T_{\om_{n-k-2j}\ldots\om_{1}}(\GJ),
\end{eqnarray*}
where we use $H(\GJ)=\GJ$.

On the other hand, we remind that $H$ acts as a contraction on trees.
Hence $T_{\be_{n'-j}\ldots \be_{1}}(\GJ)\in B(\GJ,2^{-m})$ yields
$$B(\GJ,2^{-m-1})=B\left(H(\GJ), \frac12 2^{-m}\right)\ni T_{\om_{n-k-2j}\ldots\om_{1}}(\GJ).
$$
Note that $k+2j$ is lower or equal to $2N_{m}+N_{0}$. If we set $\al=\om_{n}\ldots \om_{n-k-2j+1}$, then
$$T_{\wt\om}(\GJ)\in T_{\al}B(\GJ,2^{-(m+1)}).$$
This prove that property holds for $m+1$, and then it holds for every $m$.
\end{proof}

Following notation from \cite[Chap. IV (1.2)]{JdeVries},  Prop. \ref{prop-minimal} shows that $\GJ$ is almost periodic. Therefore, the closure of its orbit is a minimal invariant set.


\subsection{The fixed point $\GJ$ is not periodic}

We are now in a position to prove the main result of this subsection:}

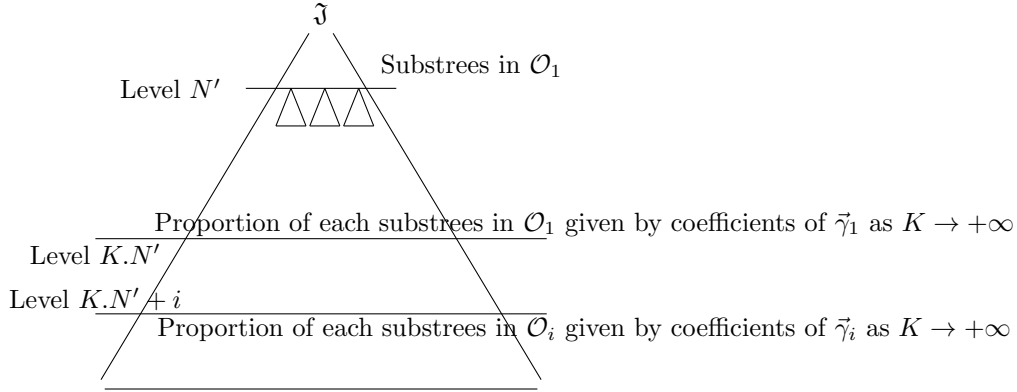
\begin{figure}[h]
\begin{tikzpicture}
\node (0) at (5,0) {$\GJ$} ;
\node (1) at (2,-5) {} ;
\node (2) at (8,-5) {} ;
\draw (0)--(1) ;
\draw (0)--(2) ;
\draw (1)--(2) ;

\draw (4,-1)--(6,-1) ;
\node at (3,-1) {Level $N'$} ;
\node at (7,-.7) {Substrees in $\CO_1$} ;

\draw (5.05,-1)--(4.85,-1.5) ;
\draw (5.05,-1)--(5.25,-1.5) ;
\draw (4.85,-1.5)--(5.25,-1.5) ;

\draw (4.6,-1)--(4.4,-1.5) ;
\draw (4.6,-1)--(4.8,-1.5) ;
\draw (4.4,-1.5)--(4.8,-1.5) ;

\draw (5.5,-1)--(5.3,-1.5) ;
\draw (5.5,-1)--(5.7,-1.5) ;
\draw (5.3,-1.5)--(5.7,-1.5) ;

\draw (2,-3)--(8,-3) ;
\node at (2,-3.2) {Level $K.N'$} ;
\node at (8.5,-2.8) {Proportion of each substrees in $\CO_1$ given by coefficients of $\vec{\ga}_{1}$ as $K\to+\8$} ;

\draw (2,-4)--(8,-4) ;
\node at (2,-3.8) {Level $K.N'+i$} ;
\node at (8.5,-4.2) {Proportion of each substrees in $\CO_i$ given by coefficients of $\vec{\ga}_{i}$ as $K\to+\8$} ;

\end{tikzpicture}
\caption{Proportion of $1$'s is bounded away from 0 as the generation increases.}
\label{fig-propsubtree}
\end{figure}

\begin{proposition}
\label{prop-fixednot perio}
The fixed point $\GJ$ is not periodic.
\end{proposition}
\begin{proof}
Because odd lines are alternations of 1 and 0, we cannot have $T_{\om}\GJ=\GJ$ with $|\om|$ odd, because it would send even lines on odd lines, and in each odd line the ration between the 0's and the 1's is equal to 1.

We do the proof by contradiction. Let us assume that there are finitely many trees $\GJ(1),\ldots,\GJ(N)$ with $\GJ(1)=\GJ$
such that the orbit of $\GJ$ is $\CO = \{ \GJ(1), \GJ(2), \ldots, \GJ(N)\}$.
This means that we can construct the graph whose vertices are the $\GJ(i)$'s and there is an arrow (with label $a$ or $b$) going from $\GJ(i)$ to $\GJ(j)$ if $\GJ(j)=T_{c}(\GJ(i))$ with $c=a,b$.
By definition of being the orbit, each vertex has two outgoing arrows and at least one incoming arrow.

\medskip
We denote by $M$ the matrix associated to this graph. If we have at least one of the equalities
$\GJ(j)=T_{a}\GJ(i)$ or  $\GJ(j) =T_{b}\GJ(i)$ then,
we set $M_{i,j}=1$; hence $M_{ij} = 1$ means that we can go from $\GJ(i)$ to $\GJ(j)$ by
at least one path ($T_a$ or $T_b$, perhaps both).

\medskip
By definition each row and each column of $M$ has at least one 1. By construction, $(M^{n})_{i,j}>0$ means that there exists $\om\in \FFF$ with $|\om|=n$ and $T_{\om}.\GJ(i)=\GJ(j)$.

We let the reader check that the matrix $M$  is exactly  of the form of the matrices considered in \cite[Ex. 1.9.4-]{katok-hasselblatt}. The spectral decomposition yields a partition of
$\cO = \left\{\GJ(1),\ldots,\GJ(N)\right\}$ in finitely many disjoint sets $\CO_{1},\ldots,\CO_{N'}$, with
$$T_{a}\CO_{i}\cup T_{b}\CO_i=\CO_{i+1},$$
with the convention $N'+1=1$.
Moreover, the dynamics for return in each component $\CO_{i}$
is mixing in the following sense: there exists $N''$ such that for every $n\ge N''$, for  every $i$, for every $k_{i}$ and $j_{i}$ such that $\GJ(k_{i})$ and $\GJ(j_{i})$ belong to $\CO_{i}$, there exists $\om\in\FFF$ with $|\om|=nN'$ and $T_{\om}\GJ(k_{i})=\GJ(j_{i})$.

\medskip
In other words,  the transition matrices for the return in each $\CO_i$ are irreducible.
We can apply the Perron-Frobenius theorem which yields for each $\CO_i$ an eigenvector with positive entries
(if we only consider the indexes $k_{i}$ such that $\GJ(k_{i})\in \CO_i$) associated to the eigenvalue
$\la^{N'}$, where $\la$ is the spectral radius\footnote{ and thus the unique single dominating eigenvalue.}
for $M$. We denote these eigenvectors  by $\vec{\ga}_{i}$.

\medskip
The $\vec{\ga}_{i}$'s yield an interesting property for the lines of $\GJ$. We remind that $\GJ=\GJ(1)$ and thus belongs to one of
the $\CO_i$'s, say $\CO_{i_{0}}$. The dynamics of the transition matrix $M$ represent the elements of the orbit of $\GJ$ that can be reached applying the $\FFF$-action. Hence, each line of value $KN'$ in $\GJ$ is only composed by roots of the elements of
$\CO_{i_{0}}$ (and they do not appear  elsewhere, since the sets $\CO_{i}$ are disjoint).
More generally, each line $KN'+i$ in $\GJ $ is only composed by roots of elements of $\CO_{(i_{0}+i\mod N')}$.

Furthermore, the entries of the vector $(1,0,\ldots, 0).M^{KN'+i}$  give exactly how many
$T_{\om_{1}\ldots \om_{KN'+i}}\GJ$ belongs to each $\GJ(k_{i})\in\CO_{(i_{0}+i\mod N')}$. In addition,  and this is a consequence of the Perron-Frobenius theorem, their relative proportion is asymptotically given by  the ratio between the entries of the eigenvectors $\vec\ga_{i_{0}+i}$ as $K$ goes to $+\8$. We refer to Fig. \ref{fig-propsubtree} for a graphic representation of that fact.
Hence, the ratio between the number of 1's and the number of 0's must be uniformly
(in $i$ and $K$) bounded away from 0 and away from 1. This is in contradiction with Corollary \ref{coro-cocorico}.
\end{proof} 

\begin{corollary}
\label{coro-JprimenondansJ}
There is no $\om\in \FFF$ such that $T_{\om}(\GJ)=\GJ'$.
\end{corollary}
\begin{proof}
The proof is done by contradiction. Let us assume that $T_{\om}(\GJ)=\GJ'$. Note that $\om\neq e$, as $\GJ$ and $\GJ'$ have different root.
Then, $T_{\om}(\GJ')=\GJ'$. Set $\om=\om_{1}\ldots \om_{n}$.

If we consider the concatenation of $k$ times $\om$, Proposition \ref{prop-minimal} shows that there exists $m\le n$ such that
$T_{\om_{1}\ldots\om_{m}}(\GJ')=\GJ$, hence $T_{\om\om_{1}\ldots\om_{m}}(\GJ)=\GJ$. This is a contradiction.
\end{proof}

\section{Self-similarities, types,  unsubstreetutions and rigidity}\label{sec-rigidity}

\subsection{Self-similarities in $\GJ$}

\begin{lemma}
\label{lem-source-brothers}
Let $\om a$ and $\om b$ be two brothers with even length. Let $\om'$ and $\om''$ respectively denote the sources of $\om a$ and $\om b$. Then either $\om'=\om''$ or they are brothers.
\end{lemma}
\begin{proof}
Let $\otimes$ be the color for the grandfather of $\om a$. We see the configuration

\begin{forest}
[{$\otimes$}[{site $\om$}[{site $\om a$}][{site $\om b$}]][{}]]
\end{forest}
or
\begin{forest}
[{$\otimes$}[{}][{site $\om$}[{site $\om a$}][{site $\om b$}]]]
\end{forest}

In any case, $\otimes$ lays at even level and we can consider its source. At the source one sees
\begin{forest}
[{$\otimes$}[{$\GA$}][{$\GB$}]]
\end{forest}
which yields at level of the grandfather
\begin{forest}
[{$\otimes$}[{1}[{$H(\GB)$}][{$H(\GB)$}]][{0}[{$H(\GA)$}][{$H(\GB)$}]]]
\end{forest}

If $\om$ is the $a$-follower of $\otimes$, then $\om'=\om''$. If $\om$ is the $b$-follower, then $\om'$ and $\om''$ are brothers.
\end{proof}

We remind that odd lines of $\GJ$ are composed by strings of $10$. We remind equalities $\chi(1)=1$, $\chi(0)=0$, $\chi(10)=0010$ and
\begin{eqnarray*}
\chi^{2}(10)&=&\chi(0010)\\
&=&\chi(10)\chi(10)\chi(00)\chi(10)\\
&=&0010\ 0010\ 0000\ 0010.
\end{eqnarray*}
Clearly, $\chi^{n}(0000)$ is a chain of 0's.

Hence, by Proposition \ref{prop-descrifixed1}, even lines in $\GJ$ are composed by blocks $0000$ or $0010$.

More precisely,
 \begin{lemma}
\label{lem-blockalbe}
Any line in $\GJ$ at level in $4\N+2$ is only composed with blocks 0010.
Any group of 16 digits in $\GJ$ at level in $4\N$ with the same $4$-ancestor is either composed only by $0$'s or is equal to
$$0010\ 0010\ 0000\ 0010.$$
\end{lemma}

\begin{lemma}
\label{lem-block0000}
Let $\om$ be a site such that all its grandchildren in $\GJ$ have root equal to $0$. Then,
$$T_{aa \wt\om}(\GJ)=T_{ba \wt\om}(\GJ)=T_{ab \wt\om}(\GJ)=T_{bb \wt\om}(\GJ).$$
\end{lemma}
\begin{proof}
Because all sites $\om aa$ $\om ab$, $\om ba$ and $\om bb$ have color $0$ it means that they lay at even level. Hence, $|\om|$ is even and we can consider its source. At the source we see
\begin{forest}
[{$\otimes$}[{$\GA$}][{$\GB$}]]
\end{forest}
which yields at site $\om$
\begin{forest}
[{$\otimes$}[{1}[{$H(\GB)$}][{$H(\GB)$}]][{0}[{$H(\GA)$}][{$H(\GB)$}]]]
\end{forest}.

This immediately yields $T_{aa \wt\om}(\GJ)=T_{ba \wt\om}(\GJ)=T_{ab \wt\om}(\GJ)=T_{bb \wt\om}(\GJ).$

Now, let us set $|\om ba|=2^{u}(2n+1)$. We can consider the unsubstreetuted trees for $H(\GB)$ and $H(\GA)$ (respectively at site $\om ba$ and $\om bb$) $u$-times. These two resulting trees are  at level $2n+1$.  By Lemma \ref{lem-source-brothers} they are either equal or ``brother'' subtrees. Furthermore their roots are equal to $0$. But at odd level, one only sees strings of $10$. This yields that the substrees are equal, hence $H(\GA)=H(\GB)$.
\end{proof}

\subsection{Types and unsubstreetution}

{ We remind $\CO(\GJ)=\disp \{T_{\om}(\GJ),\ \om\in\FFF\}$ and  $X=\ol {\CO(\GJ)}$.}

As $H(\GJ)=\GJ$, Equality \eqref{eq-renorm} shows that the set $X$ is $H$-invariant. A natural question is to inquire about in which $H^{u}(X)$ a given element belongs. If $\GA=H^{u}(\GB)$, then we say that $\GB$ is the $u$-times unsubstreetuted tree of $\GA$.
Unsubstreetution has been defined for element in $\CO(\GJ)$ in Subsection \ref{ss-unsub}, via the source function (see Lemma \ref{lem-source}). The goal of this section is to extend this notion to general elements in $X$ and not only those in the orbit of $\GJ$.

\bigskip
{ A that stage it is not yet proved that $\GJ'$ belongs to $X$. The proof is done in Corollary \ref{coro-jprimeinx} in Subsection \ref{subsec-rigi}. Nevertheless we point out that it is not an error to consider the case $\GA\in X$ and $\GA\neq\GJ'$. Furthermore, Corollary \ref{coro-jprimeinx} comes from Lemma \ref{lem-0000donne0010}. The Lemma could be stated earlier but only for $\GJ$. As we need it for any element of $X$ it appeared easier to postpone the proof of $\GJ'\in X$.}

\subsubsection{Type $2^{u}$ for elements in $X$}

\begin{definition}
\label{def-2u}
We say that $\GA\in X$ is of type $2^{u}$, with $u\in \N$ if there is a sequence $\om_{k}$ with $|\om_{k}|=2^{u}(2n_{k}+1)$ such that
$$\GA=\lim_{\kinf}T_{\om_{k}}(\GJ).$$
An element of type 1 (\ie $2^{0}$) is also said to be of odd type. An element which is not of odd type is said to be of even type.

The set of odd (resp. even) trees will be denoted by $X_{od}$ (resp. $X_{ev}$).
\end{definition}
Note that any $T_{\om}(\GJ)$ with $|\om|=2^{u}(2n+1)$ is of $2^{u}$-type. It is {\it a priori} not clear that the type is well defined for any element in $X$. The purpose of this part is to prove that it holds.

\begin{lemma}
\label{lem-defparitetype}
Let $\GA$ be in $X:=\overline{\{T_{\om}(\GJ),\ \om\in \FFF\}}$. Let $(\om_{k})$ be a sequence in $\FFF$ such that $|\om_{k}|\to_{k\to\inf}+\8$  and $\GA=\lim_{k\to\8}T_{\om_{k}}(\GJ)$.
Then only the two following cases may happen:
\begin{enumerate}
\item $(|\om_{k}|)$ eventually is even,
\item  $(|\om_{k}|)$ eventually is odd.
\end{enumerate}
\end{lemma}
\begin{proof}
We remind that any odd line for $\GJ$ is composed by a string of 10, and every even line of $\GJ$ is composed by blocks 0000 or 0010.
This yields that if $|\om|$ is even and $|\om'|$ is odd,
$$d(T_{\om}(\GJ),T_{\om'}(\GJ))\ge 2^{-3}.$$
A converging sequence is also a Cauchy sequence, thus all lengths must have the same parity for sufficiently large index.
\end{proof}

A direct consequence of Lemma \ref{lem-defparitetype} is that any $\GA\in X$ either is of odd-type or of even-type.

\begin{lemma}
\label{lem-ukinfty}
Let $\GA$ be in $X$ of even-type.
Let $(\om_{k})$ be a sequence  such that $\disp\lim_{k\to\inf}T_{\om_{k}}(\GJ)=\GA$. Set
$|\om_{k}|=2^{u_{k}}(2n_{k}+1)$. If $u_{k}\to+\8$, then $\GA=\GJ$ or $\GA=\GJ'$.
\end{lemma}
\begin{proof}
Let us set $\GA^{k}:=T_{\om_{k}}(\GJ)$. The $u_{k}$-source  for $\om_{k}$ is a site $\om'_{k}$ that is at level $2n_{k}+1$. The digit is either $0$ or $1$ (for every sufficiently large $k$ only one option is possible because $H$ is marked).
Hence, $\GA^{k}$ starts with $H^{u_{k}}(0)$ or $H^{u_{k}}(1)$. To simplify the proof we assume that the source is 0.

The assumption $u_{k}\to+\8$ yields that every integer $m$, and for every sufficiently large $k$, $A^{k}$ starts as $H^{m}(0)$. This implies that $\GA$ starts as $H^{m}(0)$, and this holds for every $m$. Hence, $\GA=\GJ$.

The same argument holds if the source is 1. In that case $\GA=\GJ'$.
\end{proof}

\begin{proposition}
\label{prop-even2u}
Let $\GA\in X$ be of even type. If $\GA\neq \GJ,\GJ'$, then it is of $2^{u}$-type for some integer $u\ge 1$ and $u$ is unique.
\end{proposition}
\begin{proof}
Let $u$ and $k$ be positive integers. Let $\om$ and $\om'$ in $\FFF$ such that $|\om|=2^{u}(2n+1)$ and $|\om'|=2^{u+k}(2m+1)$.
Set $\GB:=T_{\om}(\GJ)$ and $\GB':=T_{\om'}(\GJ)$.

Line $2^{u}$ in $\GB$ is at level $2^{u+1}(n+1)=:2^{u+a}(2l+1)$ in $\GJ$, with $a\ge1$.  By Proposition \ref{prop-descrifixed1}  this line is a concatenation of blocks $\chi^{u+a}(10)$.
Now, line $2^{u}$ in $\CB'$ is at level $2^{u}(2^{k}(2m+1)+1)$. Hence, it is a concatenation of blocks $\chi^{u}(10)$. Therefore, Corollary \ref{coro-cocorico} yields
\begin{equation}
\label{eq-2utypeloin}
d(\GB,\GB')\ge 2^{-u}.
\end{equation}

Now, consider $(\om_{k})$ such that $\GA=\lim_{k\to\8}T_{\om_{k}}(\GJ)$ with $|\om_{k}|=2^{u_{k}}(2n_{k}+1)$.
Because we have assumed $\GA\neq \GJ,\GJ'$, Lemma \ref{lem-ukinfty} shows that $(u_{k})$ is bounded. Hence, it only takes finitely many values. We may consider the smallest one which is taken infinitely many times. We denote it by $u$.

The sequence  $(T_{\om_{k}}(\GJ))$ is converging thus it is a Cauchy sequence. By \eqref{eq-2utypeloin} we must have $u_{k}\le u$ for sufficiently big $k$. This shows that the sequence $(u_{k})$ is eventually constant.

Furthermore, Corollary \ref{coro-cocorico} yields  that $u$ is unique because trees in $\GJ$ with root at lines $2^{u}(2n+1)$ and trees with root  at lines $2^{u+t}(2m+1)$ are uniformly (with respect to $n$, $m$ and position of the sites) far away from each others.
\end{proof}

\begin{notation}
In the following we shall say that $\GJ$ and $\GJ'$ are of $2^{\8}$-type.
\end{notation}

\begin{remark}
Any element of $X$ is of.a fixed type $2^{u}$ with $0\le u\le +\8$. This also holds for element in the orbit of $\GJ$. However, for elements in $\CO(\GJ)$ we also shall specify the type by the line: $T_{\om}(\GJ)$ is of the type  $2^{u}$ if $|\om|=2^{u}(2n+1)$.
\end{remark}

Now, we give two results about sequences. The next lemma is a kind of converse result from Lemma \ref{lem-ukinfty}.
\begin{lemma}
\label{lem-J=limAk}
Let $\GA^{k}:=T_{\om_{k}}(\GJ)$ be such that $|\om_{k}|=2^{u_{k}}(2n_{k}+1)$ and $\GJ=\lim_{\kinf}\GA^{k}$.
 Then $\lim_{\kinf}u_{k}=+\8$.

 The same holds if $\GJ'=\lim_{\kinf}\GA^{k}$.
\end{lemma}
\begin{proof}
Assume, by contradiction, that this is not the case. Up to a sub-sequence, we can assume that $u_{k}=u\in \N$ for every $k$. We can also assume that for every $k$, $\GA^{k}_{e}=0$ because $\lim_{\kinf}\GA^{k}=\GJ$.

Set $\GA^{k}=:H^{u}(\GB^{k})$. Each $\GB^{k}$ is an odd-type tree. They all have root equal to 0 because $H$ is marked and all $\GA^{k}$ have root equal to 0. Hence, taking the source we see a configuration
$$\GB^{k} \text{ starts as }\heptatree{0}{1,0}{0}{1}{0}{1}{0}.$$
Because $H$ is marked, $H^{u}(\GB^{k})$ differs with $\GJ$, either at line $2^{u}+1$ or at line $2^{u+1}+1$. This means that $\lim_{\kinf}\GA^{k}\neq \GJ$, which is a contradiction.

\medskip
If now we assume  $\GJ=\lim_{\kinf}\GA^{k}$.  We reemploy the same notations and mimic the proof of first point. In that case we have
$$\GB^{k} \text{ starts as }\heptatree{1}{0}{0}{1}{0}{1}{0},$$
because a1  at odd level have its two followers equal to 0.
Hence $H^{u}(\GB^{k})$ differs with $\GJ'$ at line $2^{u}+1$.
\end{proof}

\begin{proposition}
\label{prop-typeeventually}
Let $\GA$ be in $X$ and $\GA^{k}$ of type $2^{u_{k}}$ in $\CO(\GJ)$ such that $\GA=\lim_{\kinf}\GA^{k}$. Then either $u_{k}\to+\8$ or $(u_{k})$ is eventually stationary.
\end{proposition}
\begin{proof}
This is an immediate consequence of inequality \eqref{eq-2utypeloin} and the fact that a converging sequence is also a Cauchy sequence.
\end{proof}

\bigskip
In view to be able to use  Proposition \ref{prop-typeeventually}, we explain how to detect the type $2^{u}$ for any even tree.
The case $u=+\8$ is extremely simple because we must see either $\GJ$ or $\GJ'$. Hence only the finite case is relevant.

If $k$ is an integer, we denote by $v_{2}(k)$ its dyadic valuation. This means
$$v_{2}(k)=n\iff k=2^{n}(2m+1),\ m\in \N.$$
We let the reader check the following result:

\begin{lemma}
\label{lem-dyadickkpime}
Let $k, k'$ be two positive integers.
\begin{enumerate}
\item If $k'\ge k$, then for any $m$, $v_{2}(2^{k}(2m+1)+2^{k'+1})=k$.
\item If $k'=k-1$, then for any $m$, $v_{2}(2^{k}(2m+1)+2^{k'+1})\ge k+1$.
\item If $k'\le k-2$, then for any $m$, $v_{2}(2^{k}(2m+1)+2^{k'+1})= k'+1<k$.
\end{enumerate}
\end{lemma}

\begin{lemma}
\label{lem-detect-u}
Let $\GA\in X$ of $2^{u}$-type with $u<+\8$. Then, $u$ is the unique integer $k$ such that any line $2^{k+1}p$, $p\ge 1$ in $\GA$ is composed by blocks $\chi^{k}(10)$.
\end{lemma}
\begin{proof}

Let $\GA^{k}$ be a sequence of trees in $\GJ$ which converges to $\GA$. We may assume that all the roots appear at levels $2^{u}(2n_{k}+1)$. Hence, any line numbered $2^{u}2n$ in $\GA^{k}$ corresponds to a piece of line
$2^{u}(2(n_{k}+n)+1)$ in $\GJ$. By Proposition \ref{prop-descrifixed1}, such a line (in $\GJ$) is a concatenation of $\chi^{u}(10)$. This cannot be detected at line 0 in $\GA^{k}$ but it is detectable from line $2^{u+1}$ where we already have a full block $\chi^{u}(10)$.

This yields that  for any line in $\GA^{k}$ multiple of $2^{u+1}$, the proportion of 1's with respect to the 0's is fixed and is equal to the proportion of 1's in $\chi^{u}(10)$. This passes to the limit, thus this also holds in the tree $\GA$.

This shows that the set of integers $l$ such that any line $2^{l+1}n$ ($n\in\N^{*}$) is a concatenation of blocks $\chi^{l}(10)$ is non-empty and contains $u$. We now prove that $u$ is the unique integer with this property.

$\bullet$ Indeed, let us pick $v'<u$. Choose any $v$ such that $v'+1+v>u$.
Any line $2^{v'}2n$ in $\GA^{k}$ corresponds to a piece of line $2^{u}(2n_{k}+1)+2^{v'}2n=2^{v'+1}(n+2^{u-v'-1}(2n_{k}+1))$ in $\GJ$.  Now, choose a positive $n$ of the form $n=2^{v}(2m+1)-2^{u-v'-1}(2n_{k}+1)$ (with $m$ sufficiently big). Then the line $2^{v'}2n$ in $\GA^{k}$ corresponds to a piece of line $2^{v'+1+v}(2m+1)$ in $\GJ$. It is a concatenation of blocks $\chi^{v'+1+v}(10)$, but it may be truncated. As $m$ can increases as wanted, we can fix it such that the line $2^{v'}2n$ in $\GA^{k}$ is a non truncated concatenation of blocks $\chi^{v'+1+v}(10)$. Now, we remind that Corollary \ref{coro-cocorico2} states that a concatenation of blocks $\chi^{v'+1+v}(10)$ cannot be equal to a concatenation of blocks $\chi^{u}(10)$.

Because $v$ can be chosen arbitrarily, this shows that the proportion of 1's in lines $2^{v'}2n$ does not converge because it has several different accumulation points.

$\bullet$ Now for $v'>u$ and for any $n$,  a line $2^{v'}2n$ in $\GA$ is also a line $2^{u}2.2^{v'-u}n$ hence is a concatenation of blocks $\chi^{u}(10)$.

 This finishes to prove that $u$ is unique integer such that  all lines $2^{l}2n$ in $\GA$ are concatenations of blocks $\chi^{l}(10)$.

\end{proof}

\subsubsection{Unsubstreetution}

\begin{proposition}
\label{prop-unsubtree}
Let $u\ge 0$ be an integer.
Let $\GA\in X$ be of type $2^{u}$. Then there exists a unique $\GB\in X_{od}$ such that $\GA=H^{u}(\GB)$.
The tree $\GB$ is called the $u$-times unsubstreetuted tree of $\GA$, and the map $\GA\mapsto \GB$ is called the $u$-unsubstreetution.
\end{proposition}
\begin{proof}
The result is trivial if $u=0$ (\ie odd-type tree). We only do the proof for $u\ge 1$ (\ie even-type tree).

We consider a sequence $(\GA^{k})$ in $\CO(\GJ)$ converging to $\GA$. By Proposition \ref{prop-typeeventually} can assume that all the $\GA^{k}$ are of the form
$$\GA^{k}=H^{u}\left(\tritree{$\otimes$}{$\GC^{k}$}{$\GD^{k}$}\right),$$
with $\otimes=0,1$. The tree \tritree{$\otimes$}{$\GC^{k}$}{$\GD^{k}$} is of odd-type.
With these notations, we prove that $(\GC^{k})$  and $(\GD^{k})$ do converge in $X$.

The tree $H^{u}\left(\tritree{$\otimes$}{$\GC^{k}$}{$\GD^{k}$}\right)$ is of the form $H^{u}(0)$, and then at each slot at line $2^{u}+1$, we connect, either the tree $H^{u}(\GC^{k})$ or the tree $H^{u}(\GD^{k})$. If we say that $\GA^{k}$ coincides with $\GA$ for the first $n_{k}$-lines, the convergence $\GA^{k}\to \GA$ yields  $n_{k}\to+\8$ as $k\to+\8$.
Hence, $H^{u}(\GC^{k})$ and $H^{u}(\GD^{k})$
 coincide with the right
$T_{\om}(\GA)$ with $|\om|=2^{u}$ for at least $n_{k}-2^{u}$ lines.

Let us fix $N$. If we assume that for any $k\ge k_{0}$, $n_{k}>N+2^{u}$, then for any $k$ and $k'$ greater than $k_{0}$, $\GC^{k}$ and $\GC^{k'}$ on the one hand, and $\GD^{k}$  and $\GD^{k'}$ on the other hand do coincide for at least $\lfloor N2^{-u}-1\rfloor$ digits. In other words, the sequences $(\GC^{k})$  and $(\GD^{k})$ are Cauchy sequences, hence converge.
\end{proof}

We emphasize that the proof of Prop. \ref{prop-unsubtree} also yields:
\begin{lemma}
\label{lem-cvunsub}
If $\GA^{k}=H^{u}(\GB^{k})$. Then $(\GB^{k})$ converges to $\GB$ as $k\to+\8$ if and only if $(\GA^{k})$ converges to $\GA$ and $\GA=H^{u}(\GB)$.
\end{lemma}


\subsection{Rigidity}\label{subsec-rigi}

We remind that $X$ stands for the closure of the orbit of $\GJ$:
$$X:=\ol{\{T_{\om}(\GJ),\ \om\in\FFF\}}.$$

If $\GA=\tritree{$\oplus$}{$\GA'$}{$\GA''$}$ is an element in $X$, that there are actually some relations between the three quantities $\oplus$, $\GA'$ and $\GA''$. This is what we name \emph{rigidity}. Indeed, it turns out that knowing $\GA'$ or $\GA''$ actually yields knowledge of the other one. This also may lay down some specific value for $\oplus$.

\subsubsection{Blocks at even levels}

 \begin{lemma}
\label{lem-block0010}
If we see  in $X$ a tree starting as  \heptatree{$\oplus$}{$\ominus$}{$\otimes$}{$0$}{$0$}{$1$}{$0$}, then :
\begin{enumerate}
\item The tree is of even-type
\item It is equal to \heptatree{$\oplus$}{$1$}{$0$}{$\GA$}{$\GA$}{$\GB$}{$\GA$}, with $\GA_{e}=0$ and $\GB_{e}=1$.
 \end{enumerate}
\end{lemma}

\begin{proof}
Due to Lemma \ref{lem-defparitetype} and Proposition \ref{prop-typeeventually} it suffices to prove the statement for configurations in $\GJ$.
hence we consider that we see the configuration \heptatree{$\oplus$}{$\ominus$}{$\otimes$}{$0$}{$0$}{$1$}{$0$} in $\GJ$.

We denote by $\om$ the site where the root of the configuration we study stands. At line $|\om|+2$ we see a block 0010, hence line $|\om|+2$ is even. Therefore $|\om|$ is even.

Because $\om$ lays at even level, we can consider its source.  At the source we see
\begin{forest}
[{$\otimes$}[{$\GA'$}][{$\GB'$}]]
\end{forest}
which yields at site $\om$
\begin{forest}
[{$\otimes$}[{1}[{$H(\GB')$}][{$H(\GB')$}]][{0}[{$H(\GA')$}][{$H(\GB')$}]]]
\end{forest}.

This immediately yields
$T_{aa \wt\om}(\GJ)=T_{ba \wt\om }(\GJ)=T_{bb \wt\om }(\GJ)=:\GA$ and $\GB:=H(\GA')$. By hypothesis, $\GA_{e}=0$ and $\GB_{e}=1$.
\end{proof}

\begin{corollary}[Configuration 00 at even line]
\label{coro-block00}
Assume that we see in $X$ a tree \begin{forest}
[{$\oplus$}[{$\GA$}][{$\GB$}]]\end{forest}, with $\GA_{e}=\GB_{e}=0$. Then, $\oplus$ stands at odd level and $\GA=\GB$.
\end{corollary}
\begin{proof}
Again, Lemma \ref{lem-defparitetype} and Proposition \ref{prop-typeeventually}
show that it is sufficient to prove the statement for configurations  \begin{forest}
[{$\oplus$}[{$\GA$}][{$\GB$}]]\end{forest} in $\GJ$.

At the level of the roots for $\GA$ and $\GB$ we see the block 00. This shows that this level is even. By Lemma \ref{lem-blockalbe} we only see blocks 0000 or 0010 at even levels. Hence, we must see either the configuration
\begin{forest}[{}
[{$\oplus$}[{$\GA$}][{$\GB$}]]
[{}[{}][{}]]
]\end{forest}
or the reverse one
\begin{forest}[{}
[{}[{}][{}]]
[{$\oplus$}[{$\GA$}][{$\GB$}]]
]\end{forest}.
Then, Lemmas \ref{lem-block0000} and \ref{lem-block0010} show that $\GA=\GB$.
\end{proof}

\begin{lemma}[0000 yields 0010 at $2^{u}$-type, $u\ge 2$]
\label{lem-0000donne0010}
If a tree in $X$ starts as
\heptatree{$\otimes$}{$\ominus$}{$\oslash$}{$0$}{$0$}{$0$}{$0$}, then:
\begin{enumerate}
\item it is of type $2=2^{1}$.
\item we actually see \begin{forest}
[{$\otimes$}[{$1$}[{$\GA$}][{$\GA$}]][{$0$}[{$\GA$}][{$\GA$}]]]
\end{forest} with $\GA_{e}=0$.
\item $\GA$ is of $2^{u}$-type with $u\ge 2$.
\end{enumerate}

Furthermore, we also see a configuration \begin{forest}
[{$\odot$}[{$1$}[{$\GA$}][{$\GA$}]][{$0$}[{$\GB$}][{$\GA$}]]]
\end{forest} with $\GB_{e}=1$ in $X$ of the same $2^{u}$-type

\end{lemma}
\begin{proof}
Again, by Proposition \ref{prop-typeeventually} we may only prove that the Lemma holds for configurations in $\GJ$.

Then, first three statements are a direct consequence of Lemme \ref{lem-block0000}, plus the fact that a block 0000 can only appear at level in $4\N$. Note that the grand father for $\GA_{e}$ must be at an even line which is not a multiple of $4$.

Let us assume that the line where the $\GA$'s have their root is $2^{u}(2n+1)$ with $u\ge 2$.
If we consider $u$-times the source we are at level $2n+1$. We must see
\begin{figure}[h]
\begin{tikzpicture}
\Tree [.\node (level0) {$\oplus$} ; [.\node {$1$} ;
 \node {$\GE$}  ;  \node{$\GE'$}; ]
[.\node (0) {$0$}; \node (0') {$\GD$} ; \node (1) {$\GD'$};]
 ]
\draw (0) node[xshift=2cm] {level $(2n+1)$};
\end{tikzpicture}
\end{figure}

where $H^{u}\left(\begin{forest}
[{$0$}[{$\GD$}][{$\GD'$}]]\end{forest}\right)=\GA$. Actually, we have a better description because the root $\oplus$ is at even level: we see
\begin{forest}
[{$\oplus$}[
{$1$}[{$H(\GE'')$}][{$H(\GE'')$}]]
[{$0$}[{$H(\GD'')$}][{$H(\GE'')$}]]
]
\end{forest},
hence we have $\GE=\GE'=\GD'$.

Let us set $\GA':=\begin{forest}
[{$0$}[{$\GD$}][{$\GD'$}]]\end{forest}$ and $\GB':=\begin{forest}
[{$1$}[{$\GD'$}][{$\GD'$}]]\end{forest}$. This means that we see the configuration \tritree{$\oplus$}{$\GB'$}{$\GA'$}.

The image by $H$ of the subtree at level $2n$ \begin{forest}
[{$\oplus$}[{$\GB'$}][{$\GA'$}]]\end{forest} is
\begin{figure}[h]
\begin{tikzpicture}
\Tree [.\node (level0) {$\oplus$} ; [.\node {$1$} ; \node(a){$H(\GA')$};  \node(b){$H(\GA')$};  ] [.\node (r) {$0$}; \node(c){$H(\GB')$}; \node(d){$H(\GA')$}; ] ]
\draw (d) node[xshift=2cm] {level $2(2n+1)$};
\draw[dashed] (0.15,-0.8) rectangle ++(2.25,-1.5);
\end{tikzpicture}
\end{figure}

Taking the image by $H$ of the substree in the dashed frame we see
\begin{figure}[h]
\begin{tikzpicture}
\Tree [.\node (level0) {$0$} ; [.\node {$1$} ; \node(a){$H^{2}(\GA')$};  \node(b){$H^{2}(\GA')$};  ] [.\node (r) {$0$}; \node(c){$H^{2}(\GB')$}; \node(d){$H^{2}(\GA')$}; ] ]
\draw (d) node[xshift=2cm] {level $2^{2}(2n+1)$};
\draw[dotted] (0.15,-0.8) rectangle ++(2.55,-1.5);
\end{tikzpicture}
\end{figure}

In the dotted framebox, we see the same kind of configuration that previously but at level $2^{2}(2n+1)$.
Iterating this process $u-2$ more times we finally get:
\begin{figure}[h]
\begin{tikzpicture}
\Tree [.\node (level0) {$0$} ; [.\node {$1$} ; \node(a){$H^{u}(\GA')$};  \node(b){$H^{u}(\GA')$};  ] [.\node (r) {$0$}; \node(c){$H^{u}(\GB')$}; \node(d){$H^{u}(\GA')$}; ] ]
\draw (d) node[xshift=2cm] {level $2^{u}(2n+1)$};
\end{tikzpicture}
\end{figure}

Now, remember that $H^{u}(\GA')=\GA$. Moreover $\GB'_{e}=1$ and $H$ is marked, thus $H^{u}(\GB')_{e}=1$.
\end{proof}

{
\begin{corollary}
\label{coro-jprimeinx}
The tree $\GJ'$ belongs to $X$.
\end{corollary}
\begin{proof}
Minimality comes from the fact that $\GJ$ contains many finite subtrees $H^{u}(0)$ with $u$ as big as wanted. Lemma \ref{lem-0000donne0010} applied to $\GJ$ shows that any level where $\GJ$ contains $H^{u}(0)$ with $u\ge 2$, also contains a subtree starting as $H^{u}(1)$. Doing $u\to+\8$, yields that $\GJ'=\lim_{\ninf}H^{n}(1)$ belongs to $X=\ol{\CO(\GJ)}$.
\end{proof}
}

\begin{lemma}[0010 yields 0000 at $2^{u}$-type, $u\ge 2$]
\label{lem-0010donne0000}
If we see in $\GJ$ the configuration
\begin{forest}
[{$\odot$}[{$1$}[{$\GA$}][{$\GA$}]][{$0$}[{$\GB$}][{$\GA$}]]]
\end{forest}
with $\GB_{e}=1$ in $\GJ$ and $\GA$ of $2^{u}$-type with $u\ge 2$,
then we also see the configuration
\begin{forest}
[{$\otimes$}[{$1$}[{$\GA$}][{$\GA$}]][{$0$}[{$\GA$}][{$\GA$}]]]
\end{forest}
with $\GA_{e}=0$ in $\GJ$ and at the same line.
\end{lemma}
\begin{proof}
Because $\GA$ is of even-0-type, this means  $\GA_{e}$ stands at level $2^{u}(2n+1)$ with $u\ge 2$. Considering the source we see
$\begin{forest}
[{$\odot$}[{$\GB'$}][{$\GA'$}]]\end{forest}$, with $H(\GB')=\GB$ and $H(\GA')=\GA$. Roots $\GB'_{e}=1$ and $\GA'_{e}=0$ are at even level $2^{u-1}(2n+1)$. By Lemma \ref{lem-block0010} we only see blocks of 0000 or 0010. Hence we actually see

\begin{figure}[h]
\begin{tikzpicture}
\Tree [.\node (level0) {$0$} ; [.\node {$1$} ; \node(a){$\GA'$};  \node(b){$\GA'$};  ] [.\node (ro) {$0$}; \node(c){$\GB'$}; \node(d){$\GA'$}; ] ]
\draw (d) node[xshift=2cm] {level $2^{u-1}(2n+1)$};
\draw[dashed] (-1.3,-0.8) rectangle ++(1.3,-1.5);
\draw (ro) node[xshift=2cm] {with $\odot=0$} ;
\end{tikzpicture}
\end{figure}

Taking the $H$-image of the  dashed frame we get  \begin{forest}
[{$1$}[{$1$}[{$\GA$}][{$\GA$}]][{$0$}[{$\GA$}][{$\GA$}]]]
\end{forest}
 with the $\GA$'s at level $2^{u}(2n+1)$.
\end{proof}

\subsubsection{Consequences for even trees and decidability}

\begin{proposition}
\label{prop-eventree}
Only one of the three following cases may happen for an even tree $\GA$ in $X$:
\begin{enumerate}
\item $\GA=H^{v}\left(\tritree{$\GA_{e}$}{$\GD$}{$\GD$}\right)$ with \tritree{$\GA_{e}$}{$\GD$}{$\GD$} odd,
$\GD_{e}=0$ and $\GD$ of $2^{u}$-type with $1\le v\le +\8$ and $u\ge 2$.
\item $\GA=H^{v}\left(\tritree{$1$}{$\GD$}{$\GD$}\right)$ with \tritree{$1$}{$\GD$}{$\GD$} odd, $\GD_{e}=0$ and $\GD$ of $2$-type with $1\le v\le+\8$.
\item $\GA=H^{v}\left(\tritree{$0$}{$\GE$}{$\GD$}\right)$ with \tritree{$0$}{$\GE$}{$\GD$} odd, $\GD_{e}=0$  and $\GE_{e}=1$, $\GD$ and $\GE$ of type $2^{u}$, with $u\ge 1$ and $1\le v\le +\8$.
\end{enumerate}

Furthermore, knowing $T_{a}(\GA)$ or $T_{b}(\GA)$ is sufficient to determine $v$.

\end{proposition}
\begin{remark}
\label{rem-vinfty}
We emphasize that by Lemmas \ref{lem-ukinfty} and \ref{lem-J=limAk}, in the previous statement, $v=+\8$ means $\GA=\GJ$ or $\GA=\GJ'$, depending on the value for $\GA_{e}$.
$\blacksquare$\end{remark}

\begin{proof}
We consider a sequence $\GA^{k}$ in $\CO(\GJ)$ which converges to $\GA$. For simplicity we set $\GA_{e}=\otimes$. By
Proposition \ref{prop-typeeventually} we may assume that all the $\GA^{k}$ are of the form
$$H^{k}\left(\tritree{$\otimes$}{$\GE^{k}$}{$\GD^{k}$}\right),$$
and \tritree{$\otimes$}{$\GE^{k}$}{$\GD^{k}$} is an odd-tree.

By Lemma \ref{lem-cvunsub}, convergence for $\GA^{k}$ is equivalent to the convergence for $\GE^{k}$ to $\GE$ and for $\GD^{k}$ to $\GD$.

Furthermore, the root of $\GE^{k}$ and $\GD^{k}$ stand at even level, say $2^{u_{k}}(2m_{k}+1)$. Again by Proposition \ref{prop-typeeventually} we may assume that all are of fixed type $2^{u}$ with $u\ge 1$.

At even level we see blocks 0000 or 0010. Hence, the possibilities are $\GE_{e}=0=\GD_{e}$ or $\GE_{e}=1$ and $\GD_{e}=0$.
Note that by Corollary \ref{coro-block00}, $\GE_{e}=\GD_{e}$ yields $\GE=\GD$. Furthermore, $\GE_{e}=0$ is only possible if $u\ge 2$

Hence, the possibilities are
\begin{enumerate}
\item $u\ge 2$, $\GA=H^{v}\left(\tritree{$\otimes$}{$\GD$}{$\GD$}\right)$.
\item $u=1$ $\otimes=1$, $\GA=H^{v}\left(\tritree{1}{$\GD$}{$\GD$}\right)$.
\item $u\ge 1$, $\otimes=0$ and $\GA=H^{v}\left(\tritree{0}{$\GE$}{$\GD$}\right)$ with $\GE_{e}=1$ and $\GD_{e}=0$.
\end{enumerate}

\begin{remark}
\label{rem-00et10}
We point out that in the case $u\ge 2$, $\GA=H^{v}\left(\tritree{$\otimes$}{$\GD$}{$\GD$}\right)$  Lemma \ref{lem-0010donne0000} shows that the tree $\GA':=H^{v}\left(\tritree{$\bar{\otimes}$}{$\GD$}{$\GD$}\right)$ with $\bar\otimes=1-\otimes$ also exists in $X$.
$\blacksquare$\end{remark}

Now, we explain how we can determine the value of $v$. The case $v=+\8$ is extremely simple because we must see either $\GJ$ or $\GJ'$.
Hence we assume that $\GA$ is not $\GJ$ nor $\GJ'$. Then, Lemma \ref{lem-detect-u} shows that $v$ is entirely determined by one of the branches of $\GA$: for any $u$, and for any sufficiently large $n$, the line  $n$ in $T_{a}(\GA)$ or $T_{b}(\GA)$ only contains full blocks of length $|\chi^{u}(10)|$.  Hence, $v$ is the unique integer $u$ where we see blocks $\chi^{u}(10)$ at lines in $2^{u+1}\N$ (in $\GA$).

\end{proof}

\begin{remark}
\label{rem-detectles10}
We emphasize that in the case (3), $\GA$ starts as $\disp H^{v}\left(\tritree{$0$}{$1$}{$0$}\right)$. Hence lines 0 to $2^{v+1}-1$ are the same than the ones in $\GJ$. Line $2^{v}$ is equal to $\chi^v(10)$.
If $v\ge 2$, both halves of that word contain 0's and 1's (see Lemma \ref{lem-chiades1}).

Obviously, in the cases (1) and (2), the line $2^{v}$ in $\GA$ only contains 0's.
$\blacksquare$\end{remark}

\subsubsection{Branch  with root 0 determines the brother with root 1}

\begin{lemma}[Configuration 10. $0\to1$]
\label{lem-rigi0donne1}
 Assume that we see in $X$ a tree
\begin{forest}
[{$\otimes$}[{$\GA$}][{$\GB$}]]
\end{forest}, with $\GA_{e}=1$ and $\GB_{e}=0$.
Then, $\GA$ is entirely determined by $\GB$: if $\GB$ is of $2^{u}$-type and equal to $H^{u}\left(\tritree{$0$}{$\GB'$}{$\GB''$}\right)$, then $\GA=H^{u}\left(\tritree{$1$}{$\GB''$}{$\GB''$}\right)$.

\end{lemma}
\begin{proof}
Still by Proposition \ref{prop-typeeventually} we only do the proof in the case that we see the configuration in $\GJ$.

Let $\om$ be a site where we see \begin{forest}
[{$\otimes$}[{$\GA$}][{$\GB$}]]
\end{forest}.

\medskip
$\bullet$ First, we consider the case $|\om|$ is even.
If we consider the source of $\om$, we see the configuration
\begin{forest}
[{$\otimes$}[{$\GC$}][{$\GD$}]]
\end{forest}
which yields at site $\om$
\begin{forest}
[{$\otimes$}[{1}[{$H(\GD)$}][{$H(\GD)$}]][{0}[{$H(\GC)$}][{$H(\GD)$}]]]
\end{forest}.
Hence, we get $\disp \GB=\begin{forest}
[{$0$}[{$H(\GE)$}][{$H(\GD)$}]]
\end{forest}
$  and $\disp \GA=\begin{forest}
[{$1$}[{$H(\GD)$}][{$H(\GD)$}]]
\end{forest}$.

This later tree can be obtained from $\GB$
\begin{enumerate}
\item by exchanging the color of the root,
\item replacing the left hand side subtree by a copy of the one at the right-hand side.
\end{enumerate}
This corresponds to the result of the Lemma with $u=0$.

\medskip
$\bullet$ Assume that $|\om|$ is odd. Then we set $|\om|+1=2^{u}(2n+1)$. We consider the $u$-source for $\GA$ and $\GB$. By Lemma \ref{lem-source-brothers} they are either brothers or equal. As $H$ is marked and $\GA_{e}=1\neq 0=\GB_{e}$, the $u$-sources are brothers. Let us denote them by $\GA^{u}$ and $\GB^{u}$. They are odd-type trees at level $2n+1$. Furthermore, $\GA^{u}_{e}=1$ and $\GB^{u}_{e}=0$. By Lemma \ref{lem-blockalbe}, their father is 0.

We use the part of Lemma \ref{lem-rigi0donne1} that is already proved (with father at even level). Hence $\GA^{u}$ is obtained from $\GB^{u}$ by exchanging the color of the root and replacing $T_{a}(\GB^{u})$ by $T_{b}(\GB^{u})$. Then, applying $H^{u}$ we get for $\GA$:

\begin{enumerate}
\item $\GA$ starts as $H^{u}(1)$ whereas $\GB$ starts with $H^{u}(0)$.
\item At any slot at line $2^{u}$ we glue $H^{u}(T_{b}(\GB^{u}))$.
\end{enumerate}
Employing notation from the statement of the lemma we have
$\GB=H^{u}\left(\GB^{u}\right)$, $\GA=H^{u}\left(\GA^{u}\right)$.
$\GB^{u}=\tritree{$0$}{$\GB'$}{$\GB''$}$, then $\GA^{u}=\tritree{$1$}{$\GB''$}{$\GB''$}$.
\end{proof}

\begin{remark}
\label{rem-coherent}
We emphasize that the construction of $\GA$ from $\GB$ in Lemma \ref{lem-rigi0donne1} is coherent with the construction of $\GB$ from $\GA$ in Lemma \ref{lem-0000donne0010}.
$\blacksquare$\end{remark}


\section{Preimages and proof of Theorem \ref{MainThB}}\label{sec-mainB}

\subsection{Preimages in general}
Given a tree $\GA$ in $X$, the preimages are of the form \tritree{$\otimes$}{$\GA$}{$\GB$} or \tritree{$\otimes$}{$\GB$}{$\GA$}. We remind that $\GA\in X$ means $\GA=\lim_{k\to+\8}\GA^{k}$, $\GA^{k}=T_{\om_{k}}(\GJ)$.

The question of preimages is then equivalent to determine for a given converging sequence $(\GA^{k})$, which configurations \tritree{$\otimes$}{$\GA^{k}$}{$\GB^{k}$} or \tritree{$\otimes$}{$\GB^{k}$}{$\GA^{k}$} we do see in $\GJ$.

We shall see that answers depend on the type of the configuration $\GA$: even or odd.
It also depends on the values for the roots of $\GA^{k}$ and $\GB^{k}$.

From Lemmas  \ref{lem-block0000}, \ref{lem-blockalbe} and \ref{lem-block0010}, for even type trees we only see\footnote{the root level is the framed one.}
\begin{forest}
[{0,1}
[{1}[0, draw][0, draw]]
[{0}[{0,1}, draw][0, draw]]
]
\end{forest}.
 For odd-type trees, we only see\footnote{the root level is the framed one.}
\begin{forest}
[{0,1}[1, draw][0,draw]]
\end{forest}.

\begin{enumerate}
\item In Subsection \ref{ss1-odd} we deal with odd-type trees. Proposition \ref{prop-preimimparac0} deals with odd trees with root equal to 0 and Proposition \ref{prop-preimimparac1} deals with odd trees with root equal to 1.
\item In Subsection \ref{ss-even0} we deal with even type trees with root equal to 0. Proposition \ref{prop-preim0pairrac0} deals with $2^{u}$-types with $2\le u<+\8$, Proposition \ref{prop-preim1pairrac0} deals with $2$-type, and Proposition \ref{pre-jaca} deals with the case $2^{\8}$-type ($\GJ$ and $\GJ'$).
\item In Subsection \ref{ss-even1} we deal with the even type trees and root equal to 1. Proposition \ref{prop-preim0pairrac1} deals with the $2^{u}$-types, $1\le u<+\8$ and Proposition \ref{prop-T-1Jprime} deals with the case $2^{u}$-type $u=+\8$.
\end{enumerate}
The proof of Theorem \ref{MainThB} is done in Subsection \ref{ss-proofB}


\subsection{Preimage of odd trees}\label{ss1-odd}

Let $\GA, \GB \in X$ be odd trees with $\GA_{e}=1$ and $\GB_{e}=0$. Then :
\begin{enumerate}
\item  $T^{-1}(\GA)=T_{a}^{-1}(\GA)$. \\
\item $T^{-1}(\GB)=T_{b}^{-1}(\GB)$.
\end{enumerate}
In the first case a brother is an odd tree $\GB'$ with $\GB'_{e}=0$ and such that \tritree{$\otimes$}{$\GA$}{$\GB'$} exists in $X$. In the second case, a brother is an odd tree $\GA'$ with $\GA_{e}=1$ and such that \tritree{$\otimes$}{$\GA'$}{$\GB$} exists in $X$.

\medskip

Because the root for odd trees determine if there are or not preimages by $T_{a}$ or $T_{b}$, we shall, in the subsection keep the letter $\GA$ for odd trees with root 1 and $\GB$ for odd trees with root 0.

\medskip

We emphasize that for the later case, Lemma \ref{lem-rigi0donne1} shows that $\GA'$ is entirely determined by $\GB$.

\bigskip

We now see that the same holds in the first case. First we point out that an odd tree with root equal to  1  is of the form \tritree{1}{$\GD$}{$\GD$} with  $\GD_{e}=0$, $\GD$ of $2^{u}$-type with $1\le u\le+\8$.

{ Furthermore, any brother for $\GA$ is a tree of the form \tritree{0}{}{$\GD$}. More precisely, Lemmas \ref{lem-0000donne0010} and \ref{lem-0010donne0000} show that the configurations \tritree{0}{1}{$\GD$} always exists  but the configuration \tritree{0}{0}{$\GD$} does exists if and only if $\GD$ is of $2^{u}$-type with $u\ge 2$. This is formulated more explicitly in the following lemma. }

\begin{lemma}
\label{lem-brotherA}
Let $\GA\in X$ be an odd tree with root equal to 1. Set $\GA=\tritree{$1$}{$\GD$}{$\GD$}$ with  $\GD_{e}=0$, $\GD$ of $2^{u}$-type with $1\le u\le+\8$.

Let $\GF$ be the tree obtained applying Lemma \ref{lem-rigi0donne1} to \tritree{$0$}{$\GF$}{$\GD$}.
Set $\GB:=\tritree{$0$}{$\GF$}{$\GD$}$ and $\GB':=\tritree{$0$}{$\GD$}{$\GD$}$

Then,

\begin{enumerate}
\item Whatever $u\ge 1$ is, there exists a tree \tritree{}{$\GA$}{$\GB$} in $X$.
\item If  (and only if) $u\ge 2 $, there exists a tree \tritree{}{$\GA$}{$\GB'$} in $X$.
\end{enumerate}
 \end{lemma}
\begin{proof}

As $\GA$ belongs to $X$, we see in $\GJ$ trees \heptatree{$\otimes$}{1}{0}{$\GD^{k}$}{$\GD^{k}$}{.}{.} with $\lim_{\kinf}\GD^{k}=\GD$. The trees $\GD^{k}$ stand at $2^{u}$-levels.
By Lemma  \ref{lem-0000donne0010}), we will actually see in $\GJ$ trees \heptatree{$\otimes$}{1}{0}{$\GD^{k}$}{$\GD^{k}$}{$\GD^{k}$}{$\GD^{k}$} (this holds only for $u\ge 2$) and trees \heptatree{$\odot$}{1}{0}{$\GD^{k}$}{$\GD^{k}$}{$\GF^{k}$}{$\GD^{k}$} (at the same levels) with $\GF^{k}_{e}=1$ (this holds for any $u\ge 1$).

The tree $\GF^{k}$ is obtained from $\GD^{k}$ via the procedure described in Lemma \ref{lem-rigi0donne1}.
Lemma \ref{lem-cvunsub} yields that $\GF^{k}$ converges as $\GD^{k}$ converges.  The link between $\GF^{k}$ and $\GD^{k}$ passes to the limit and then $\GF^{k}$ converges to $\GF$ (from the statement of the Lemma).

Therefore, doing $\kinf$ we get:

\begin{enumerate}
\item in any case,  a tree \tritree{$\odot$}{$\GA$}{$\GB$} in $X$,
\item if and only if $u\ge 2$, a tree \tritree{$\otimes$}{$\GA$}{$\GB'$} in $X$.
\end{enumerate}

\end{proof}

Lemma \ref{lem-brotherA} states (for a given $\GA\in X_{od}$ with root 1) what kind of trees  \tritree{$\otimes$}{$\GA$}{$\GB$} we can find in $X$. However, it does not make precise what is the value for $\otimes$. Similarly, we have seen  that Lemma \ref{lem-rigi0donne1} states that an odd tree $\GB$ has a single brother $\GA'$ in $X$. Again, it does not make precise what the father may be.

Furthermore, preimages of these trees are even trees, and thus obey the rule described in Proposition \ref{prop-eventree}. We now connect the two results and give a full description of preimages of odd trees.

\begin{proposition}
\label{prop-preimimparac1}
Let $\GA$ be an odd tree satisfying $\GA_{e}={\RL 1}$. Let $v$ be such that $\GA\in T_{a}\circ H^{v}(X_{od})$. Then,
\begin{enumerate}
\item If $v\ge 2$, then $\GA$ has a unique brother $\GB$ given by Lemma \ref{lem-brotherA} case (1). Furthermore,
only one of the following cases may appear:
\begin{enumerate}
\item $T_{a}^{-1}(\GA)=\left\{\tritree{$0$}{$\GA$}{$\GB$},\tritree{$1$}{$\GA$}{$\GB$}\right\}$, line $2^{v}-1$ in $\GA$ only contains 0's and $T^{2^{v}-1}(\GA)$ is even of $2^{v+u}$-type with $u\ge 2$.
\item $T_{a}^{-1}(\GA)=\left\{\tritree{$1$}{$\GA$}{$\GB$}\right\}$, line $2^{v}-1$ in $\GA$ only contains 0's  and $T^{2^{v}-1}(\GA)$ is even of $2^{v+1}$-type.
\item $T_{a}^{-1}(\GA)=\left\{\tritree{$0$}{$\GA$}{$\GB$}\right\}$, line $2^{v}-1$ in $\GA$ contains 0's  and 1's.
\end{enumerate}
\item If $v=1$, then $\GA$ has two brothers $\GB$ and $\GB'$ given by Lemma \ref{lem-brotherA}. Furthermore, only one of the following cases may appear:
\begin{enumerate}
\item $T_{a}^{-1}(\GA)=\left\{\tritree{$0$}{$\GA$}{$\GB'$},\tritree{$1$}{$\GA$}{$\GB'$}, \tritree{$0$}{$\GA$}{$\GB$}\right\}$, and $T_{a}(\GA)$ is of $2^{1+u}$-type with $u\ge 2$.
\item $T_{a}^{-1}(\GA)=\left\{\tritree{$1$}{$\GA$}{$\GB'$}, \tritree{$0$}{$\GA$}{$\GB$}\right\}$, and $T_{a}(\GA)$ is  of $2^{2}$-type.
\end{enumerate}
\end{enumerate}

\end{proposition}
\begin{proof}
We remind that by Proposition \ref{prop-eventree}, $\GA$ determines the type of any preimage.

$\bullet$ If $v\ge 2$ holds, then any tree \tritree{}{$\GA$}{} starts as \heptatree{$\otimes$}{1}{0}{0}{0}{1}{0}, which shows that $\GA$ has a unique brother, and it starts as \tritree{0}{1}{0}.

Using notations from Lemma \ref{lem-brotherA}, this yields that $\GB$ is the unique brother for $\GA$ (case $\GB'$ is impossible).

Line $2^{v}-1$ in $\GA$ corresponds to the line $2^{v} in $\tritree{}{$\GA$}{$\GB$}. By Remark \ref{rem-detectles10}, this line only contains 0's if and only if we are in the cases (1) or (2) from Proposition \ref{prop-eventree}, and the line contains 0's and 1's if and only if we are in the case (3) from Proposition \ref{prop-eventree}.

Furthermore, for any $\om\in \FFF$ with $|\om|=2^{v}-1$, $T_{\om}(\GA)$ is equal to $H^{v}(\GD)$ if we re-employ notations from Proposition \ref{prop-eventree}. It is thus of type $2^{u+v}$ and this determine the type of $\GD$. Hence, this differentiates case (1) and  case (2) (from Proposition \ref{prop-eventree}).

Now, still for $v\ge 2$, case (1) in Proposition \ref{prop-eventree} corresponds to (1-a) in the statement of Proposition \ref{prop-preimimparac1}, case (2) in Proposition \ref{prop-eventree} corresponds to (1-b) in the statement of Proposition \ref{prop-preimimparac1},
 and case (3) in Proposition \ref{prop-eventree} corresponds to (1-c) in the statement of Proposition \ref{prop-preimimparac1},

\bigskip
$\bullet$ Let us now assume that $v=1$ holds. Reemploying notations from Lemma \ref{lem-brotherA} we see \tritree{}{$\GA$}{$\GB$} and \tritree{}{$\GA$}{$\GB'$} in $X$. Line 2 in \tritree{}{$\GA$}{$\GB$} is equal to 0010, whereas line 2 in \tritree{}{$\GA$}{$\GB'$} is equal to 0000.
Hence, the tree \tritree{}{$\GA$}{$\GB$} corresponds to the case (3) from Proposition \ref{prop-eventree}, and \tritree{}{$\GA$}{$\GB$} corresponds to the case (1) or (2) from Proposition \ref{prop-eventree}.

If we set $\GA=\tritree{$1$}{$\GG$}{$\GG$}$. Because $v=1$,  \tritree{}{$\GA$}{} belongs to $H(X_{od})$ hence $\GG$ is of $2^{k}$-type  with $k\ge 2$. This means that if we set $\disp \tritree{}{$\GA$}{$\GB'$}=H\left(\tritree{}{$\GD$}{$\GD$}\right)$, then $\GD$ is of type $k-1$.

 If $k-1=1$, only the configuration \tritree{1}{$\GD$}{$\GD$} does exist because lines at even level but not in $4\N$ only contains block 0010.  On the contrary, if $k-1\ge 2$, then Lemmas \ref{lem-0000donne0010} and \ref{lem-0010donne0000} yield that both configurations \tritree{1}{$\GD$}{$\GD$} and \tritree{0}{$\GD$}{$\GD$} do exist.

Hence we get: If $T_{a}(\GA)$ is of type $k\ge 3$, then $\GA$ has three preimages, that are \tritree{0}{$\GA$}{$\GB'$}
, \tritree{1}{$\GA$}{$\GB'$} and \tritree{0}{$\GA$}{$\GB$}.

If $T_{a}(\GA)$ is of type $k=2$, then $\GA$ has two preimages, that are \tritree{1}{$\GA$}{$\GB'$}
and \tritree{0}{$\GA$}{$\GB$}.
\end{proof}

\begin{proposition}
\label{prop-preimimparac0}
Let $\GB\in X$ be an odd tree with $\GB_{e}=0$. Let $\GA$ be its brother (as Lemma \ref{lem-rigi0donne1}). Set $T_{b}(\GB)$ of type $2^{k}$. Let $v$ be such that \tritree{}{$\GA$}{$\GB$} is of type $2^{v}$ (as in Prop. \ref{prop-eventree}). Then, only one of the following cases may appear:
\begin{enumerate}
\item If $\GB$ starts as \tritree{$0$}{$0$}{$0$}, then
\begin{enumerate}
\item $k\ge 3$ and $v=1$  and $\disp T^{-1}_{b}(\GB)=\left\{\tritree{$0$}{$\GA$}{$\GB$},\tritree{$1$}{$\GA$}{$\GB$}\right\}$.
\item $k=2$ and $v=1$ and \tritree{$1$}{$\GA$}{$\GB$} is the unique preimage of $\GB$.
\end{enumerate}
{ \item If $\GB$ starts as \tritree{$0$}{$1$}{$0$}, then
\begin{enumerate}
\item If $k\ge 2$, then $v=1$ and $\disp T^{-1}_{b}(\GB)=\left\{\tritree{$0$}{$\GA$}{$\GB$}\right\}$.
\item If $k=1$ then $v\ge 2$ and  one of the following case occurs:
\begin{enumerate}
\item  line $2^{v}-1$ of $\GB$ only contains 0's and $T^{2^{v}-1}(\GB)$ is of $2^{u+v}$-type with $u\ge 2$. Then $\disp T^{-1}_{b}(\GB)=\left\{\tritree{$0$}{$\GA$}{$\GB$},\tritree{$1$}{$\GA$}{$\GB$}\right\}$.
\item  line $2^{v}-1$ of $\GB$ only contains 0's and $T^{2^{v}-1}(\GB)$ is of $2^{1+v}$-type. Then $\disp T^{-1}_{b}(\GB)=\left\{\tritree{$1$}{$\GA$}{$\GB$}\right\}$.
\item Line $2^{v}-1$ of $\GB$ contains 0's and 1's. Then $\disp T^{-1}_{b}(\GB)=\left\{\tritree{$0$}{$\GA$}{$\GB$}\right\}$.
\end{enumerate}
\end{enumerate}}
\end{enumerate}
\end{proposition}
\begin{proof}
Any tree \tritree{$\otimes$}{$\GA$}{$\GB$} is of the form $\disp H^{v}\left(\tritree{$\otimes$}{$\GE$}{$\GD$}\right)$ with \tritree{$\otimes$}{$\GE$}{$\GD$} odd. We remind that $v$ is uniquely determined by $\GB$ (see Prop. \ref{prop-eventree}).

$\bullet$ If $\GB$ starts\footnote{which means that it is a $\GB'$ if we re-employ notations from Lemma \ref{lem-brotherA}.} as \tritree{$0$}{$0$}{$0$}. This is possible only if $v=1$, as $H^{2}(\otimes)$ starts as \heptatree{$\otimes$}{1}{0}{0}{0}{1}{0}. Hence we get $k=u+1$. Furthermore,

-- either $u\ge 2$ (which is equivalent to $k\ge 3$) and then $\GB=\disp T_{b}H\left(\tritree{$0$}{$\GD$}{$\GD$}\right)=\disp T_{b}H\left(\tritree{$1$}{$\GD$}{$\GD$}\right)$ for some $\GD$ of type $2^{u}$ satisfying $\GD_{e}=0$. Note that by Lemmas \ref{lem-0000donne0010} and \ref{lem-0010donne0000}, both configurations \tritree{1}{$\GD$}{$\GD$} and \tritree{0}{$\GD$}{$\GD$} do exist.
 This means that $\GB$ has 2 preimages.

-- or $u=1$ (\ie $k=2$) and $\GB=\disp T_{b}H\left(\tritree{$0$}{$\GD$}{$\GD$}\right)$  for some $\GD$  of $2^{1}$-type with $\GD_{e}=0$. Here only the configuration \tritree{1}{$\GD$}{$\GD$} does exist. This means that $\GB$ has a unique preimage with root 1.

$\bullet$ If $\GB$ starts as \tritree{$0$}{$1$}{$0$}. Note that \tritree{$\otimes$}{$\GA$}{$\GB$} starts as \heptatree{$\otimes$}{1}{0}{0}{0}{1}{0}. If $T_{a}(\GB)$ is of $2^{k}$-type with $k\ge 2$, then $\otimes$ stands at an even level which is not in $4\N$. This yields $v=1$.

Hence we get $\tritree{$\otimes$}{$\GA$}{$\GB$}=H\left(\tritree{$\otimes$}{$\GE$}{$\GF$}\right)= \heptatree{$\otimes$}{$1$}{$0$}{$H(\GF)$}{$H(\GF)$}{$H(\GE)$}{$H(\GF)$}$. This yields $\GE_{e}=1$ and $\GF_{e}=0$. Now, \tritree{$\otimes$}{$1$}{$0$} with $\otimes$ at odd level is possible in $\GJ$ if and only if $\otimes=0$.

\medskip
Now, we deal with the case $\GB$ starts as \tritree{$0$}{$1$}{$0$} and $k=1$. Remind that \tritree{$\otimes$}{$\GA$}{$\GB$} starts as \heptatree{$\otimes$}{1}{0}{0}{0}{1}{0} and the last line $0010$ stands at a level in $2\N\setminus 4\N$. This yields that the root $\otimes$ stands at a level in $4\N$, hence $v\ge 2$. We remind that $v$ can be detected, as it is explained in Lemma \ref{lem-detect-u}. Furthermore,  Proposition \ref{prop-eventree} states that $v$ can be detected knowing only $\GB$.

 Hence we have:

--either we are in the cases (1) or (2) from Prop. \ref{prop-eventree} with $v\ge 2$. In that case there is a line of 0's at line $2^{v}$ of \tritree{$\otimes$}{$\GA$}{$\GB$}.

-- or we are in the case (3) from Prop. \ref{prop-eventree} (with $v\ge 1$). In that case the line $2^{v}$ of \tritree{$\otimes$}{$\GA$}{$\GB$} contains 0's and 1's.

Furthermore, to separate  cases (1) and (2), we need to check what is the type $2^{u+v}$ of any tree $T_{\om}(\GB)$ with $\om\in \FFF$ and $|\om|=2^{v}-1$.
\end{proof}

\newpage
 The table \ref{tab-preimaimpar} summarizes the possibles preimages for a given odd tree.

\begin{table}[h]
\begin{tabular}{|p{3cm}||p{3cm}|p{3cm}|p{3cm}|p{1cm}|p{2 cm}|}
  \hline
  $\swarrow$ $T_{a} $ or $T_{b}$ & $H^{v}\left(\tritree{$0,1$}{$\GD$}{$\GD$}\right)$, $\GD$ of $2^{u}$-type $u\ge 2$, $\GD_{e}=0$ & $H^{v}\left(\tritree{$1$}{$\GD$}{$\GD$}\right)$, $\GD$ of $2$-type, $\GD_{e}=0$ & $H^{v}\left(\tritree{$0$}{$\GF$}{$\GD$}\right)$, $\GD_{e}=0$, $\GF_{e}=1$ & Number of preimages& Case\\
  \hline \hline
  \multirow{6}{3cm}{$\GA$, $\GA_{e}=1$ $T_{a}(\GA)$ of $2^{k}-type$} & $v\ge 2$, brother $\GB$ &&&2&Prop. \ref{prop-preimimparac1}-1-a\\
 \cline{2-6}
  &&$v\ge 2$, brother $\GB$ &&1&Prop. \ref{prop-preimimparac1}-1-b\\
  \cline{2-6}
&&&$v\ge 2$, brother $\GB$&1& Prop. \ref{prop-preimimparac1}-1-c\\
  \cline{2-6}
& $v=1$ and brother $\GB'$ $k=u+1$&  & $v=1$, $k=u+1$ brother $\GB$&3&Prop. \ref{prop-preimimparac1}-2-a\\
\cline{2-6}
&&$v=1$ and brother $\GB'$ $k=2$ & $v=1$, $k=2$ brother $\GB$&2&Prop. \ref{prop-preimimparac1}-2-b\\
 \hline\hline
  $\GB'$,  starts as \tritree{0}{0}{0},  $T_{a}(\GB')$ of $2^{k}-type$& $k\ge 3$, $v=1$ & & &2&Prop. \ref{prop-preimimparac0}-1-a\\
 \cline{2-6}
 &&$k=2$, $v=1$&&1&Prop. \ref{prop-preimimparac0}-1-b\\
  \hline\hline
 $\GB$,  starts as \tritree{0}{1}{0}, $T_{a}(\GB)$ of $2^{k}-type$&$v\ge 2$&&&2&Prop. \ref{prop-preimimparac0}-2-b-i\\
  \cline{2-6}
  &&$v\ge 2$&&1&Prop. \ref{prop-preimimparac0}-2-b-ii\\
   \cline{2-6}
  &&&$v\ge 2$&1&Prop. \ref{prop-preimimparac0}-2-b-iii\\
   \cline{2-6}
 &&&{\RL$v=1$}&1&Prop. \ref{prop-preimimparac0}-2-a\\
  \hline
\end{tabular}
\captionof{table}{Preimages of odd trees. All cases in the same line coexist. }\label{tab-preimaimpar}-2-a
\end{table}


\subsection{Preimage for even type elements with root equal to 0}\label{ss-even0}

\begin{proposition}
\label{prop-preim0pairrac0}
Let $\GB$ be in $X$ of $2^{u}$-type with $2\le u<+\8$ and $\GB_{e}=0$. Then,
$T^{-1}(\GB)$ is the set
$$\left\{\begin{forest}
[{$0$}[{$\GB$}][{$\GB$}]]\end{forest}, \begin{forest}
[{$1$}[{$\GB$}][{$\GB$}]]\end{forest}, \begin{forest}
[{$0$}[{$\GA$}][{$\GB$}]]\end{forest}
\right\},$$
where  $\GA=H^{u}\left(\tritree{$1$}{$\GB''$}{$\GB''$}\right)$ if $\GB=H^{u}\left(\tritree{$0$}{$\GB'$}{$\GB''$}\right)$.
\end{proposition}
\begin{proof}
We consider a sequence of trees $\GB^{k}:=T_{\om_{k}}(\GJ)$ such that $\GB=\disp\lim_{\kinf}\GB^{k}$. For simplicity we assume that for every $k$, $|\om_{k}|=2^{u}(2n_{k}+1)$ holds.

Lemmas \ref{lem-0000donne0010} and \ref{lem-0010donne0000} shows that we see the configurations \begin{forest}
[{$\odot$}[{$1$}[{$\GB^{k}$}][{$\GB^{k}$}]][{$0$}[{$\GA^{k}$}][{$\GB^{k}$}]]]
\end{forest}
and \begin{forest}
[{$\otimes$}[{$1$}[{$\GB^{k}$}][{$\GB^{k}$}]][{$0$}[{$\GB^{k}$}][{$\GB^{k}$}]]]
\end{forest},
with $\GA^{k}_{e}=1$. Moreover, $\GA^{k}$ is obtained from $\GB^{k}$ following Lemma \ref{lem-0000donne0010} and as explained in  Lemma \ref{lem-rigi0donne1}.

As $\GB^{k}$ goes to $\GB$ if $\kinf$, branches at level $2^{u}$ for $\GB^{k}$ converge. Hence $\GA^{k}$ converges to some tree $\GA\in X$.
This yields
$$T^{-1}(\GB) \supset \left\{\begin{forest}
[{$0$}[{$\GB$}][{$\GB$}]]\end{forest}, \begin{forest}
[{$1$}[{$\GB$}][{$\GB$}]]\end{forest}, \begin{forest}
[{$0$}[{$\GA$}][{$\GB$}]]\end{forest}
\right\}.$$

Furthermore, configurations \tritree{$1$}{$\GB^{k}$}{$\GD^{k}$} or \tritree{$1$}{$\GD^{k}$}{$\GB^{k}$} with $\GD^{k}_{e}=0$ and $\GD^{k}\neq \GB^{k}$ are forbidden.
In the same direction, configurations  \tritree{$0$}{$\GB^{k}$}{$\GD^{k}$} with $\GD^{k}_{e}=0$ and $\GD^{k}\neq \GB^{k}$ are also forbidden. And finally, configurations \tritree{$0$}{$\GB^{k}$}{$\GD^{k}$} with $\GD^{k}_{e}=1$ are forbidden.

\end{proof}

\begin{proposition}
\label{prop-preim1pairrac0}
Let $\GB$ be in $X$  2-type and $\GB_{e}=0$. Then,
$T^{-1}(\GB)$ is the set
$$\left\{\begin{forest}
[{$1$}[{$\GB$}][{$\GB$}]]\end{forest}, \begin{forest}
[{$0$}[{$\GA$}][{$\GB$}]]\end{forest}
\right\},$$
where $\GA=H\left(\tritree{$1$}{$\GB''$}{$\GB''$}\right)$ if $\GB=H\left(\tritree{$0$}{$\GB'$}{$\GB''$}\right)$.
\end{proposition}
\begin{proof}
The proof is almost the same than for  Prop. \ref{prop-preim0pairrac0}. The unique difference is that, by Proposition \ref{prop-descrifixed1}, the configuration
\begin{forest}
[{$\otimes$}[{$1$}[{$\GB^{k}$}][{$\GB^{k}$}]][{$0$}[{$\GB^{k}$}][{$\GB^{k}$}]]]
\end{forest}
does not exist.

\end{proof}

\begin{proposition}
\label{pre-jaca}
$$T^{-1}(\GJ)= \left\{\begin{forest}
[{$0$}[{$\GJ$}][{$\GJ$}]]\end{forest}, \begin{forest}
[{$1$}[{$\GJ$}][{$\GJ$}]]\end{forest}, \begin{forest}
[{$0$}[{$\GJ'$}][{$\GJ$}]]\end{forest}
\right\}.$$
\end{proposition}
\begin{proof}
The fixed point $\GJ$ is not periodic. Hence, there is no $\om\neq e$ such that $T_{\om}(\GJ)=\GJ$. This yields that there is no $\om\in \FFF$ with $\om\neq e$ such that $T_{\om}(\GJ)$ belongs to $T^{-1}(\GJ)$.

An element $\GC$ in $T^{-1}(\GJ)$ is thus obtained as $\lim_{\kinf}\GC^{k}$ with $\GC^{k}$ of the form $T_{\om_{k}}(\GJ)$.
If, say $T_{a}(\GC)=\GJ$, then $\lim_{\kinf} T_{{\RL a\om_{k}}}(\GC^{k})=\GJ$.

Conversely, if $\GB^{k}$ is of the form $\GB^{k}=T_{\al_{k}}(\GJ)$ and $\lim_{\kinf}\GB^{k}=\GJ$, if $\GC^{k}$ is such that, say $T_{a}\GC^{k}=\GB^{k}$, then any accumulation point $\GC$ for the $\GC^{k}$'s satisfies $T_{a}(\GC)=\GJ$.

\medskip
To compute $T^{-1}(\GJ)$, we thus adapt the proof of Prop. \ref{prop-preim0pairrac0}. We consider a sequence  of terms $\GB^{k}:=T_{\om_{k}}(\GJ)$ such that $\lim_{\kinf}\GB^{k}=\GJ$.  We compute what $T^{-1}(\GB^{k})$ is and consider all the accumulation points for these sequences of trees if $k$ goes to $+\8$.

Note that $\GB^{k}_{e}=0$ (at least for any sufficiently large $k$). By Lemma \ref{lem-J=limAk}, if we set $|\om_{k}|=2^{u_{k}}(2n_{k}+1)$, then $u_{k}\to+\8$ as $k\to+\8$. Hence, for any sufficiently large $k$, $u_{k}\ge 2$ and we can apply Prop. \ref{prop-preim0pairrac0} to each $\GB^{k}$. For simplicity we assume that the above conditions hold for every $k$.

Then, for each fixed $k$, $T^{-1}(\GB^{k})$ is the set whose elements are \tritree{0}{$\GB^{k}$}{$\GB^{k}$}, \tritree{1}{$\GB^{k}$}{$\GB^{k}$} and \tritree{0}{$\GA^{k}$}{$\GB^{k}$}, where $\GA^{k}$ is given in Lemma \ref{lem-rigi0donne1}.

The first 2 elements converge respectively to \tritree{0}{$\GJ$}{$\GJ$} and \tritree{1}{$\GJ$}{$\GJ$} as $k$ goes to $+\8$.
Now, we remind that $\GA^{k}$ starts as $H^{u_{k}}(1)$ and $u_{k}\to+\8$. Furthermore, in the proof of Prop. \ref{prop-Jprime} we have seen that $H^{n}(0)$ and $H^{n}(1)$ differ only from the root. Therefore, $\lim_{\kinf}\GA^{k}=\GJ'$.
This shows that \tritree{0}{$\GJ'$}{$\GJ$} is also in $X$ and in $T^{-1}(\GJ)$.
\end{proof}

\subsection{Preimages of even-type elements with root equal to 1}\label{ss-even1}

We start with a description of even-type elements with root equal to 1.

\begin{lemma}
\label{lem-descrievenroot1}
Let $u\ge 1$ be an integer.
$2^{u}$-type elements with root equal to 1 in $X$ are elements of the form $H^{u}\left(\tritree{$1$}{$\GC$}{$\GC$}\right)$ with $\GC=H^{v}\left(\tritree{$0$}{$\GC'$}{$\GC''$}\right)$, $1\le v\le +\8$.
\end{lemma}
\begin{remark}
\label{rem-vinftyaeven}
The case $v=+\8$ means $\GC=\GJ$.
$\blacksquare$\end{remark}
\begin{proof}
Because $\GA$ is of even-type and $\GA_{e}=1$ and $H$ is marked, Prop. \ref{prop-unsubtree} yields that we can write $\GA$ under the form $\GA=H^{u}(\GB)$ with $\GB_{e}=1$, and $\GB$ is of odd-type. Any 1 at odd line is followed only by 0's and Corollary \ref{coro-block00}  yields $$\GB=\tritree{$1$}{$\GC$}{$\GC$},$$
with $\GC_{e}=0$.

The tree $\GB$ is of odd-type, thus $\GC$ is of even-type, say $2^{v}$-type with $1\le v\le +\8$. Hence, $\GC$ can be written under the form \tritree{$0$}{$\GC'$}{$\GC''$}, with the convention that $v=+\8$ means $\GC=\GJ$.

\end{proof}

\begin{proposition}
\label{prop-preim0pairrac1}
Let $\GA$ be in $X$ of $2^{u}$-type with $1\le u<+\8$ and $\GA_{e}=1$.
Set $\GA=H^{u}\left(\tritree{$1$}{$\GC$}{$\GC$}\right)$ with $\GC=H^{v}\left(\tritree{$0$}{$\GC'$}{$\GC''$}\right)$, $1\le v\le +\8$. Then,
\begin{enumerate}
\item If $v=1$, then
$$T^{-1}(\GA)=T_{a}^{-1}(\GA)=\left\{\begin{forest}
[{$0$}[{$\GA$}][{$\GB$}]]\end{forest}\right\}$$
with $\GB=H^{u}\left(\tritree{$0$}{$\GD$}{$\GC$}\right)$ and $\GD=H^{v}\left(\tritree{$1$}{$\GC''$}{$\GC''$}\right)$.
\item If $v>1$, then
$$T^{-1}(\GA)=T_{a}^{-1}(\GA)
=\left\{\begin{forest}
[{$0$}[{$\GA$}][{$\GB$}]]\end{forest},\begin{forest}
[{$0$}[{$\GA$}][{$\GB'$}]]\end{forest}\right\},$$
with $\GB=H^{u}\left(\tritree{$0$}{$\GD$}{$\GC$}\right)$, $\GD=H^{v}\left(\tritree{$1$}{$\GC''$}{$\GC''$}\right)$ and
$\GB'=H^{u}\left(\tritree{$0$}{$\GC$}{$\GC$}\right)$.
\end{enumerate}
\end{proposition}
\begin{remark}
\label{rem-sivinfty}
We remind that if $v=+\8$ then we replace $H^{v}\left(\tritree{$0$}{}{}\right)$ by $\GJ$ and $H^{v}\left(\tritree{$1$}{}{}\right)$ by $\GJ'$.
Note that $\GB'$ differs from $\GA$ only with the root.
$\blacksquare$\end{remark}

\begin{proof}
We consider a sequence of terms $\GA^{k}:=T_{\om_{k}}(\GJ)$, with $\GA^{k}\to \GA$ as $k\to+\8$. With the notations of the Lemma, we may assume that every $\GA^{k}$ is of $2^{u}$-type and we set $\GA^{k}=H^{u}\left(\tritree{$1$}{$\GC^{k}$}{$\GC^{k}$}\right)$. By Lemma \ref{lem-cvunsub} $\GC^{k}\to\GC$. Hence, we may write $\GC^{k}=H^{v}\left(\tritree{$0$}{$\GG^{k}$}{$\GH^{k}$}\right)$. This makes sense if $v\neq +\8$. If $v=+\8$, then we set $\GC^{k}=H^{v_{k}}\left(\tritree{$0$}{$\GG^{k}$}{$\GH^{k}$}\right)$ with $v_{k}\to+\8$.

The root of $\GA^{k}$ (equal to 1) is at an even line. We remind that even lines are composed by blocks 0000 or 0010. This shows that $T_{b}^{-1}(\GA^{k})=\emptyset$.

By Lemma \ref{lem-block0010} we see a  sequence of configurations
 \begin{forest}
[{$\otimes$}[{$1$}[{$\GB^{k}$}][{$\GB^{k}$}]][{$0$}[{$\GA^{k}$}][{$\GB^{k}$}]]]
\end{forest},
with $\GB^{k}_{e}=0$.

Note that this is the unique configuration where $\GA^{k}$ appears, and this yields that $T^{-1}(\GA)=T^{-1}_{a}$ and that any preimage of $\GA$ has its root equal to 0.

\medskip
We now see how the $\GB^{k}$'s do depend on the $\GA^{k}$'s.

We remind equality $\GA^{k}=H^{u}\left(\tritree{$1$}{$\GC^{k}$}{$\GC^{k}$}\right)$. Hence, we consider the $u$-times unsubstreetuted trees for $\GA^{k}$ and $\GB^{k}$. They have different roots since $H$ is marked. By  Lemma \ref{lem-source-brothers} they are brothers.

Hence we see a configuration
\heptatree{$\otimes$}{1}{0}{$\GC^{k}$}{$\GC^{k}$}{$\GM^{k}$}{$\GN^{k}$}.

We first deal with the case $v<+\8$.
The root of $\GC^{k}$ stands at even level $2^{v}(2m_{k}+1)$ and its value is 0. Because at even line we only see blocks of 0010 or blocks of 0000, Lemma \ref{lem-block0010} yields $\GN^{k}=\GC^{k}$.

Lemmas \ref{lem-0000donne0010} and \ref{lem-0010donne0000} shows that:
\begin{enumerate}
\item Either $v=1$ and necessarily $\GM^{k}_{e}=1$.
\item Or $v\ge 2$ and we see at the same level both configurations \heptatree{$\otimes$}{1}{0}{$\GC^{k}$}{$\GC^{k}$}{$\GM^{k}$}{$\GC^{k}$} for some $\GM^{k}$  with $\GM^{k}_{e}=1$ and \heptatree{$\otimes$}{1}{0}{$\GC^{k}$}{$\GC^{k}$}{$\GC^{k}$}{$\GC^{k}$}.
\end{enumerate}

For both cases where we see some subtree $\GM^{k}$, it remains to determine how we obtain the subtree $\GM^{k}$ from the subtree $\GC^{k}$. For that we use Lemma \ref{lem-rigi0donne1}. We have $\GC^{k}=H^{v}\left(\tritree{$0$}{$\GG^{k}$}{$\GH^{k}$}\right)$, thus
$$\GM^{k}=H^{v}\left(\tritree{$1$}{$\GH^{k}$}{$\GH^{k}$}\right).$$
Because $\GC^{k}$ converges to $\GC=\tritree{$0$}{$\GC'$}{$\GC''$}$, Lemma \ref{lem-cvunsub} shows that $\GG^{k}$ converges to (say) $\GC'$ and $\GH^{k}$ converges to $\GC''$.  Furthermore, the same lemma shows that $\GM^{k}$ converges to $H^{v}\left(\tritree{$1$}{$\GC''$}{$\GC''$}\right)$.

Going back to $\GB^{k}$ we finally get:
\begin{enumerate}
\item If $v=1$, there is a unique preimage for $\GA$, it is \tritree{0}{$\GA$}{$\GB$} with $\GB=H^{u}\left(\tritree{$0$}{$\GD$}{$\GC$}\right)$ where $\GD=H^{v}\left(\tritree{1}{$\GC''$}{$\GC''$}\right)$.
\item If $2\le v<+\8$, then $\GA$ has two preimages. The same than for the case $v=1$ plus \tritree{0}{$\GA$}{$\GB'$} with $\GB'=H^{u}\left(\tritree{$0$}{$\GC$}{$\GC$}\right)$
\end{enumerate}

\bigskip
Now, we deal with the case $v_{k}\to+\8$. In that case $\GA=H^{u}\left(\tritree{$1$}{$\GJ$}{$\GJ$}\right)$. The same work as for the case $v>1$ can be done, except that we finally consider $v\to+\8$.
Hence we get two preimage for $\GA$. The tree  \tritree{0}{$\GA$}{$\GB$} with $\GB=H^{u}\left(\tritree{$0$}{$\GJ'$}{$\GJ$}\right)$ and the tree
\tritree{0}{$\GA$}{$\GB'$} with $\GB'=H^{u}\left(\tritree{$0$}{$\GJ$}{$\GJ$}\right)$.

\end{proof}

\begin{proposition}
\label{prop-T-1Jprime}
$T^{-1}(\GJ')=T_{a}^{-1}(\GJ')=\left\{\tritree{$0$}{$\GJ'$}{$\GJ$}\right\}$.
\end{proposition}
\begin{proof}
By corollary \ref{coro-JprimenondansJ}, there is no $\om$ such that $T_{\om}(\GJ)=\GJ'$. Therefore, preimages of $\GJ'$ are the accumulation points of preimages of $\GA^{k}:=T_{\om_{k}}(\GJ)$ such that $\lim_{\kinf}\GA^{k}=\GJ'$.

Let us consider such a sequence, and set $\om_{k}=2^{u_{k}}(2m_{k}+1)$. By Lemma \ref{lem-J=limAk}, $\lim_{\kinf}u_{k}=+\8$.

If we reemploy notations from the beginning of the proof of Prop. \ref{prop-preim0pairrac1}, we will see in $\GJ$ the configuration \begin{forest}
[{$\odot$}[{$1$}[{$\GB^{k}$}][{$\GB^{k}$}]][{$0$}[{$\GA^{k}$}][{$\GB^{k}$}]]]
\end{forest}.
By Lemma \ref{lem-source-brothers}, $\GA^{k}$ starts has $H^{u_{k}}(1)$ whereas $\GB^{k}$ starts as $H^{u_{k}}(0)$. Then, by Lemma \ref{lem-ukinfty}, $\lim_{\kinf}\GB^{k}=\GJ$. This finishes the proof.
\end{proof}

\subsection{Proof of Theorem \ref{MainThB}}\label{ss-proofB}
The proof is immediate. It follows from Propositions \ref{prop-preimimparac0}, \ref{prop-preimimparac1}, \ref{prop-preim0pairrac0}, \ref{prop-preim1pairrac0},  \ref{prop-preim0pairrac1} and \ref{prop-T-1Jprime} that for any $\GA\in X$,
 $$p(1,\GA)=\#T^{-1}(\GA)\cap X\le 3.$$
 Hence, it yields $p(n,\GA)\le 3^{n}$.

\section{Some more results}\label{sec-Appendix}

Here we present some examples from grammars others than BBAB.

\subsection{Thue-Morse case}

Define the function $H$ in such a way that

$$
\disp 0\mapsto
\begin{forest}
[$0$[$1$][$1$]]ullet
\end{forest}  \;\; \text{and} \;\; \disp 1\mapsto
\begin{forest}
[$1$[$0$][$0$]]
\end{forest},
$$

Then using the recurrence it is easy to see that  fixed points for  $H$ are ``grammar-free'': each one of the fixed points $\afz$ and $\afu$ has each generation with only one symbol. There is thus a natural projection $\pi$ from $\afz$ and $\afu$ to $\{0,1\}^{\NN}$, commuting the dynamics. Moreover, $\pi\circ H(0)=01$ and $\pi\circ H(1)=10$. This means that $\pi\circ H$ is the classical Thue-Morse sequence.

Hence, $\overline{\{T_{\om}(\afz),\ \om\in\FFF\}}=\overline{\{T_{\om}(\afu),\ \om\in\FFF\}}$ is a minimal dynamical system, although it does not really exploit the tree structure.

\subsection{Proof of Theorem  \ref{th-nonminimal}}

We consider the substreetution given by
$$
\disp 0\mapsto  \begin{forest}
[$0$[$0$][$1$]]ullet
\end{forest}  \;\; \text{and} \;\; \disp 1\mapsto \begin{forest}
[$1$[$1$][$0$]]
\end{forest},
$$
equipped with the grammar ABBA.
It is marked and therefore admits two fixed points, $\afz$ and $\afu$, respectively with root $0$ and $1$.

The renormalization equation is
\begin{equation} \label{renormeqABBA}
   T_a T_a H = H T_a = T_b  T_b H \qquad \text{and} \qquad  T_a T_b H = T_b T_a H = H T_b
\end{equation}
and then the source is given by
$$\begin{cases}
s(aa)=s(bb)=a,\\
s(ab)=s(ba)=b.\\
\end{cases}$$
More generally we have
\begin{equation}
\label{recABBAeven}
s(p_1 p_2 \ldots p_{2n}) =  s(p_1 p_2) \ldots s(p_{2n-1} p_{2n}) .
\end{equation}

 Using the identification $a = 0$ and $b=1$ we can write simply $s(xy) = (x+y) \mod(2)$.

The images by $H^{2}$ of roots $0$ and $1$ are

\begin{forest}
[$0$[$0$, [$0$[$0$][$1$]][$1$[$1$][$0$]]][$1$[$1$[$1$][$0$]][$0$[$0$][$1$]]]]
\end{forest}
and
\begin{forest}
[$1$[$1$, [$1$[$1$][$0$]][$0$[$0$][$1$]]][$0$[$0$[$0$][$1$]][$1$[$1$][$0$]]]]
\end{forest}

In the following, $\GD$ stands for $\afz$ or $\afu$.
 We note the following equalities
$$
  H(\GD)_e  = \GD_e,  H(\GD)_a = \GD_e, H(\GD)_b = 1-\GD_e = \bar{\GD}_e
$$
(where $\bar{x} = 1- x$).
We remark that {\bf identifying $a$ with $0$ and $b$ with $1$} then the last two  expressions
above can be rewritten as
$$
  H(\GD)_{\om} = (\GD_e + \om) \mod(2)   \qquad \text{where $\om = 0 , 1$}.
$$
By induction, we get  for the odd generations
\begin{equation} \label{recABBAodd}
  \GD_{p_1 p_2 \ldots p_{2n} a} =  \GD_{p_1 p_2 \ldots p_{2n} }  \quad \text{and} \quad
   \GD_{p_1 p_2 \ldots p_{2n} b} =  \bar{\GD}_{p_1 p_2 \ldots p_{2n} }.
\end{equation}
These expressions can be rewritten as
$$
\GD_{p_1 p_2 \ldots p_{2n} \om} =  ( \GD_{p_1 p_2 \ldots p_{2n} } + \om)  \mod(2) ,
$$
where we still make the identification $a=0$ and $b=1$.

\begin{lemma}
Writing $\GD = \afi$, $i=0,1$, then the following holds:
$$
  \forall n\ge 1,\  \GD_{p_1 p_2 \ldots p_n} = (\GD_e + p_1 + p_2 + \cdots + p_n) \mod(2)  =  (i + p_1 + p_2 + \cdots + p_n) \mod(2).
$$
\end{lemma}
\begin{proof}
The proof is by induction.
For $n =1$ it is true, from the definition of $H$.

Assume that it is valid for $n$; then $\GD_{p_1 \ldots p_k} =  (\GD_e + p_1 + p_2 + \cdots + p_k) \mod(2)$
for any $k=1, \ldots, n$.

Now we need to show that the expression is valid for $n+1$.
If $n+1$ is odd the we can use (\ref{recABBAodd}) to obtain
$$
 \GD_{p_1 p_2 \ldots p_n p_{n+1}} = ( \GD_{p_1 \ldots p_n} + p_{n+1}) \mod(2) = ( \GD_e + p_1 + \cdots + p_n + p_{n+1} ) \mod(2)
$$
as wished.

If, on the other hand, $n+1$ is even then $n+1 = 2k$ and so we have, from (\ref{recABBAeven}), that
$$
  \GD_{p_1 p_2 \ldots p_{n+1}} = \GD_{p_1 p_2 \ldots p_{2k-1} p_{2k}} = \GD_{f(p_1 p_2) f(p_3 p_4) \cdots f(p_{2k-1}p_{2k}) } =
$$
$$
   (\GD_e + f(p_1 p_2) + \ldots + f(p_{2k-1} p_{2k}) ) \mod(2)  =
$$
$$
   = (\GD_e + p_1 + p_2 + p_3 + p_4 + \ldots + p_{2k-1} + p_{2k} ) \mod(2) =
$$
$$
  = (\GD_e + p_1 + p_2 + \ldots + p_{n} + p_{n+1}) \mod(2)
$$
as claimed, concluding the proof.
\end{proof}

The expression above allows us to write a generation from the previous one, just using that $0$ gives rise to $01$
(on the subsequent generation) and $1$ gives rise to $10$ (also on the subsequent generation).
This gives us an easy recursive way to write the two fixed points of $H$ in this case.   The fixed
points $\afz$ and $\afu$ (whose roots are, respectively, $0$ and $1$) then can be obtained writing the generation
$n+1$ from the generation $n$. In this case we then  have
 $T_a (\afz) = T_b (\afu) = \afz$ and $T_b (\afz) = T_a (\afu) = \afu$.

An easy consequence is that for every $n\ge 1$, $\afz_{ba^{n}}=1$. This prevents $\afz$ to be almost periodic. By \cite[Chap. IV (1.2)]{JdeVries},  $\overline{\{T_{\om}(\afz),\ \om\in \FFF\}}$ is not minimal. By symmetry, the same holds for $\overline{\{T_{\om}(\afu),\ \om\in \FFF\}}$.

\subsection{Proof of Theorem \ref{th-perio-nomeas}}
We give here an example of periodic tree whose orbit does not support an invariant measure.
The tree is not obtained by a subtreetution but via another way: each line is obtained from the previous one by applying some map. The choice of the map depends on the parity of the line.

We set
$$\begin{cases}
\CT_{1}(0) = 01\text{ and } \CT_{1}(1) = 10,\\
\CT_{2}(01)=0001\text{ and }\CT_{2}(10)=1110.\\
\end{cases}$$

We call $\afz$ and $\afu$ the two colored trees obtained by the following process with roots $0$ and $1$:
an even line $2k$ is obtained applying $\CT_2$ to the line $2k-1$ (gluing digit by pair); a line
$2k+1$ is obtained by applying $\CT_1$ to each digit of the line $2k$.

The fixed points starts as
\begin{forest}
[{\hskip-0.8cm$\afz=0$}[$0$[$0$[$0$][$1$]][$0$[$0$][$1$]]][$1$[$0$[$0$][$1$]][$1$[$0$][$1$]]]]
\end{forest}
and
\begin{forest}
[{\hskip-0.8cm$\afu=1$}[$1$[$1$[$1$][$0$]][$1$[$1$][$0$]]][$0$[$1$[$1$][$0$]][$0$[$1$][$0$]]]]
\end{forest}

As for each even line we apply $\CT_1$, which depends only on one digit, then the subtree only depends on
that site. Hence at each even line, at each site with color $0$ there is a $\afz$ and at each site of color $1$
there is a $\afu$. This means that we have
$\afz =$ \heptatree{0}{0}{1}{$\afz$}{$\afz$}{$\afz$}{$\afu$} and
 $\afu =$ \heptatree{1}{1}{0}{$\afu$}{$\afu$}{$\afu$}{$\afz$}

Hence
$$
 T_a (\afz) = \tritree{$0$}{$\afz$}{$\afz$} =:\GB, \; T_b(\afz) = \tritree{$1$}{$\afz$}{$\afu$} =: \GC, \;
  T_a(\afu) = \tritree{$1$}{$\afu$}{$\afu$} =: \GD, \; T_b(\afu) = \tritree{$0$}{$\afu$}{$\afz$} =: \GE
$$
and
$$
  T_a(\GB) = T_b(\GB) = \afz, T_a(\GC) = \afz, T_b(\afz)=\afu, T_a(\GD)=T_b(\GD) = \afu, T_a(\GE)=\afu, T_b(\GE)=\afz
$$
Then the orbit of $\afz$ is $\CO = \{\afz, \afu, \GB, \GC, \GD, \GE  \}$ and $T_a(\CO) \cup T_b(\CO) = \CO$.

\begin{figure}[h]\label{fig-arbnomeasure}
\begin{tikzpicture}

\node[draw,circle] (b)at(-2,0) {$\GB$};

\node[draw,circle] (0)at(0,0) {$\afz$} ;
\node[draw,circle] (1)at(4,0) {$\afu$} ;
\node[draw,circle] (c)at(2,2) {$\GC$} ;
\node[draw,circle] (e)at(2,-2) {$\GE$} ;

\node[draw,circle] (d)at(6,0) {$\GD$};

\draw[->,>=latex] (0) to[bend left=20] node[midway]{$T_{b}$} (c);
\draw[->,>=latex] (c) to[bend left=20] node[midway]{$T_{a}$} (0);

\draw[->,>=latex] (c) to[bend left=20] node[midway]{$T_{b}$} (1);

\draw[->,>=latex] (0) to[bend left=20] node[midway]{$T_{a}$} (b);
\draw[->,>=latex] (b) to[bend left=20] node[midway]{$T_a,T_{b}$} (0);

\draw[->,>=latex] (e) to[bend left=20] node[midway]{$T_{b}$} (0);
\draw[->,>=latex] (e) to[bend left=20] node[midway]{$T_a$} (1);

\draw[->,>=latex] (1) to[bend left=20] node[midway]{$T_b$} (e);
\draw[->,>=latex] (1) to[bend left=20] node[midway]{$T_a$} (d);

\draw[->,>=latex] (d) to[bend left=20] node[midway]{$T_a, T_b$} (1);



\end{tikzpicture}
\caption{The graph for $\FFF$-action on $\CO$}
\end{figure}

We remind the reader that a measure $\mu$ is said to be an invariant probability if for every Borel set $A$ we have
$\mu(A) = \mu(T_a^{-1}(A)) = \mu(T_b^{-1}(A))$. In particular, $\mu(\GC) = \mu(T_a^{-1}(\GC)) = \mu(\emptyset) = 0$;
similarly we get $\mu(\GB) = \mu(\GD) = \mu(\GE) = 0$. Finally,
$\mu(\afz) = \mu(T_a^{-1}(\afz)) = \mu(\GB \cup \GC) = 0$ and   $\mu(\afu) = \mu(T_a^{-1}(\afu)) = \mu(\GD \cup \GE) = 0$,
showing that we do not have an invariant probability on the periodic orbit $\CO$.

\subsection{Colored quasi-periodic tilings for hyperbolic disk}

In \cite{BH} it is proved that there does not exist a primitive substitutive tiling of the hyperbolic plane $\mathbb{H}^{2}$. Tilings studied there are geometrical ones. We show here how our colored trees may generate colored tilings for the hyperbolic plane. The construction is done in the hyperbolic disk.

The following example of  ``ping-pong'' in the Hyperbolic disk has kindly been given to us by M. Peign\'e.

We consider in $\mathbb{D}^{2}$ the isometry $h_{1}$ which maps $-\frac12$ to $\frac12$ and let $\pm1$ fixed\footnote{We use the canonical embedding of $\mathbb{D}^{2}$ into $\mathbb C$}.
We also consider the isometry $h_{2}$ which maps $-\frac{i}2$ to $\frac{i}2$ and let $\pm i$ fixed.

The Dirichlet domain (see \cite{dalbo}) is obtained by considering the intersection of the half-planes
$$\{z\in \mathbb{D}^{2}, d(z,0)\le d(z,\ga 0)\},$$
with $\ga=h_{i}^{\pm1}$, $i=1,2$.

\begin{figure}[htbp]
\begin{center}
\includegraphics[scale=0.5]{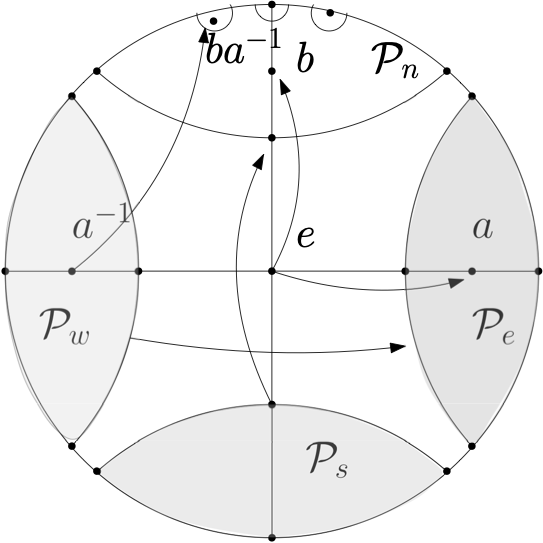}
\caption{Colored tiling from ping-pong in $\mathbb{D}$}
\label{fig-pingpong.eps}
\end{center}
\end{figure}

We consider the half-planes $\{z\in \mathbb{D}^{2}, d(z,0)\le d(z,\ga 0)\},$ with $\ga=h_{i}^{\pm1}$, $i=1,2$. They are denoted by $\CP_{w}$, $\CP_{e}$, $\CP_{s}$ and $\CP_{n}$ as on Figure \ref{fig-pingpong.eps}. We set $\CP_{0}:=\mathbb{D}^{2}\setminus\cup_{i=w,e,s,n}\CP_{i}$. Hence $\CP_{0}$ is the Dirichlet domain we deal with.

Furthermore, we identify $0\sim e$ $a\sim h_{1}$ and $b\sim h_{2}$.

\medskip
Now, for any given $\GA\in X$ we construct a colored tiling for $\mathbb{D}^{2}$ in the following way: For each $\om\in \FFF$, the fundamental domain that contains $\om$ in $\mathbb{D}^{2}$ (up to identification) has the color of the site $\GA_{\om}$.
Each half-plane $\om.b(\CP_{w})$ (resp. $\om.a(\CP_{s})$) has the same color than $\GA_{\om.b}$ (resp. $\GA_{\om_.a}$).

This tiling only colors images $\om.\CP_{0}$ with $\om\in \FFF$. We can also decide to color all different $\om^{-1}.\CP_{0}$ via the same principle. It can also use another tree $\GB\in X$, up to the condition that the color for $\CP_{0}$ has to be fixed (either $\GA_{e}$ or $\GB_{e}$).

\subsection{Picture of Jacaranda tree}

The next picture has kindly been  computed by F. Flouvat.
We point out that the tree structure increases exponentially fast the number of sites and this is a brake to compute  the fixed point on a large number of generations. Graphic representation is also a challenge. The next picture  used the algorithm introduced in \cite{palgrave}.

\begin{figure}[htbp]
\begin{center}
\includegraphics[scale=0.7]{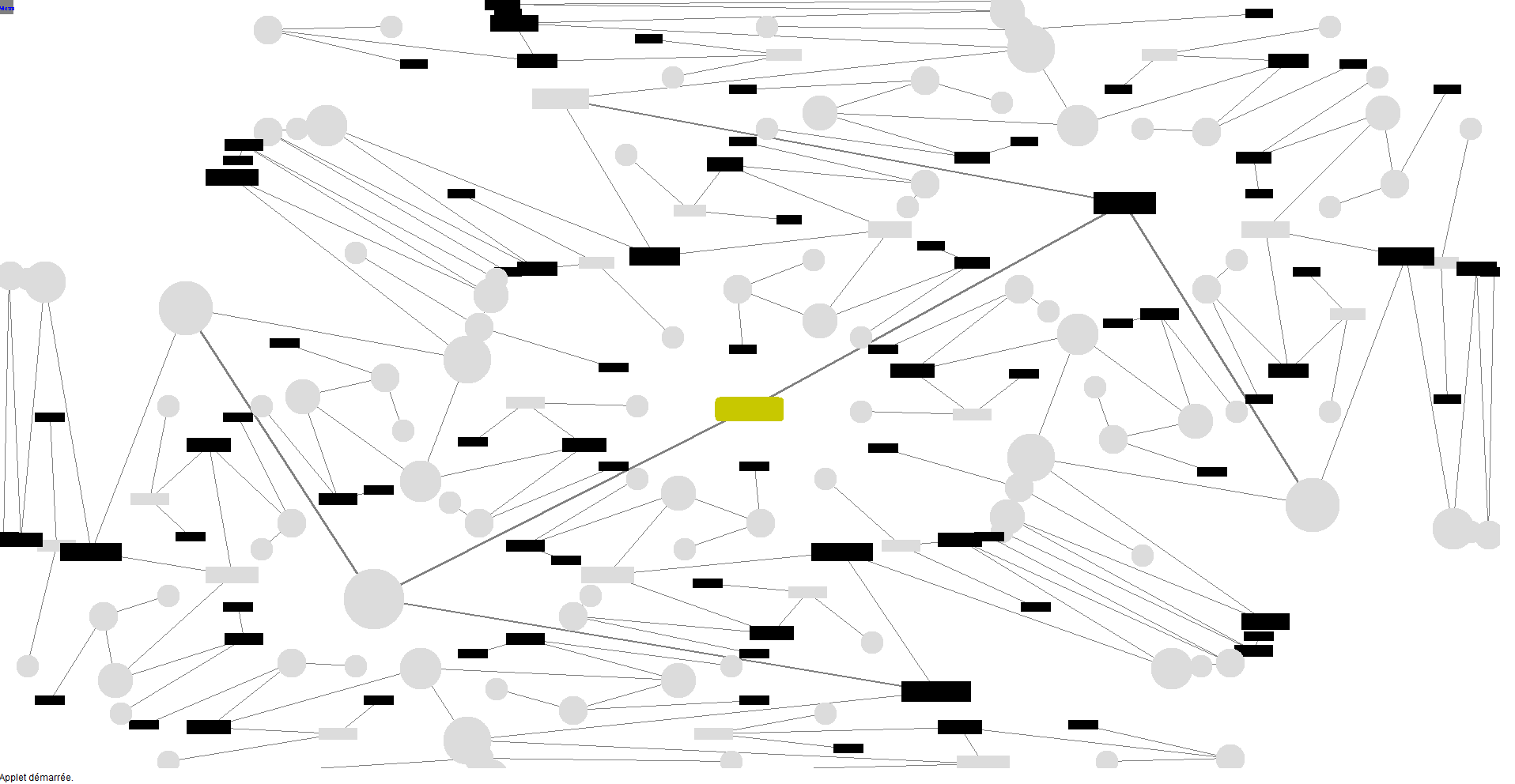}
\caption{The Jacaranda tree: black=1, grey =0, $a$-follower= rectangle, $b$-follower= disk, root=yellow}
\end{center}
\end{figure}



\begin{thebibliography}{27}

\bibitem{BBCF}  J. Berstel, L. Boasson, O. Carton and I. Fagnot, \emph{A first investigation of Sturmian Trees},
 Feb 2007, Aachen (Aix la Chapelle), Germany. pp.73-84,
10.1007/978-3-540-70918-3 7 . hal-00150178

\bibitem{BH}  N. B\'edaride and A. Hilion, \emph{Geometric realizations of 2-dimensional substitutive tilings},  Q. J. Math. {\bf 64} (2013), no. 4, 955–-979.



\bibitem{BLL} A. Baraviera, R. Leplaideur and A. Lopes,
\emph{The potential point of view for Renormalization},  Stoch. \& Dynam. {\bf 12} (2012), Issue : 4


\bibitem{BL1} H. Bruin and R. Leplaideur,
\emph{Renormalization, thermodynamic formalism  for quasi-crystals in subshifts},
 Commun. Math. Phys. {\bf 231} (2013), pp. 209-–247

\bibitem{BL2} H. Bruin and R. Leplaideur,
\emph{Renormalization, freezing phase transition and Fibonacci quasicrystals},
 Annales Scientifiques de l’ENS {\bf 48} (2015), fascicule 3, pp. 739--763

\bibitem{BP} I. Benjamini and Y. Peres, 
\emph{Markov chains indexed by trees}, Ann. Probability {\bf 22} (1994), no. 1, 219-243 


\bibitem{dalbo}
{Dal'Bo, F.},
   \emph{Geodesic and horocyclic trajectories}, {Universitext}, {Springer-Verlag London, Ltd., London; EDP Sciences, Les Ulis}, {2011}.
   	
\bibitem{Damm} W. Damm, \emph{The IO- and OI-Hierarchies}
	Theoretical Computer Science 20: 95-207 (1982)

\bibitem{Dek} F. M. Dekking, \emph{On the Thue-Morse Measure}, Acta Univ. Carolinae {\bf 33} (1992), no. 2, 35--40

\bibitem{Fog} P. Fogg, \emph{ Substitutions in dynamics, arithmetics and combinatorics}, Lecture Notes in Mathematics {\bf 1794} (2002)

\bibitem{FW} H. Furstenberg and B. Weiss, \emph{Markov Processes and Ramsey theory for trees}, Combinatorics,
Probability and Computing  {\bf 12} (2003), pp. 547--563


\bibitem{JdeVries} J. de Vries, \emph{Elements of topological dynamics}, Springer (1993)

\bibitem{katok-hasselblatt} A. Katok  and B. Hasselblatt, \emph{Introduction to the
modern theory of dynamical systems}, Cambridge University Press (2010)

\bibitem{KLLS} D. H. Kim, B. Lee, S. Lim and D. Sim, \emph{Quasi-Sturmian colorings on regular trees}, pre-print arXiv

\bibitem{palgrave} Q. V. Nguyen and M. L. Huang, \emph{Space optimized tree: a connection+enclosure approach for
the visualization of large hierarchies}, Information Visualization {\bf 2} (2003), 3--15


\bibitem{Petersen-Salama-20}
{Petersen, Karl and Salama, Ibrahim},
     \emph{Entropy on regular trees},
 {Discrete Contin. Dyn. Syst.},
     {40},
      {2020},
  { \bf no.7},
   {4453--4477},


\bibitem{Przytycki-pressure} F. Przytycki. \emph{Conical limit set and Poincaré exponent for iterations of rational functions}
Trans. Amer. Math. Soc.
351 (1999), {\bf no. 5}, 2081–2099

\bibitem{PRLS}
F. Przytycki J. Rivera-Letelier and
S. Smirnov,  \emph{Equality of pressures for rational functions}
Ergodic Theory Dynam. Systems,
24, (2004)
{\bf no.3}, 891--914


\bibitem{Que} M. Queffelec, \emph{ Substitution dynamical systems-spectral analysis}, Lecture Notes in Mathematics {\bf 1294} (2010)

\bibitem{Rozikov}  U. A. Rozikov, \emph{Gibbs measures on Cayley Trees}, World Scientific (2013)

\bibitem{Spitzer} F. Spitzer, \emph{Markov random fields on an infinite tree}, Ann. Probability {\bf 3} (1975), no. 5, 387-398

\end{thebibliography}
\end{document}